\documentclass[11pt, a4paper]{article}
\usepackage{verbatim}
\usepackage{amsmath,latexsym}
\usepackage{amssymb,wasysym}
\usepackage[latin1]{inputenc}
\usepackage{soul}%define \st para tachare
\usepackage{slashbox}
\usepackage{ marvosym }
\usepackage{array}
\usepackage{tabularx}
\usepackage{fancyhdr}
\usepackage{psfrag}
\usepackage{xspace}
\usepackage{pstricks}
\usepackage{pst-node,xcolor,graphicx}
\usepackage{pst-tree}
\usepackage{graphicx}
\usepackage{subfigure}
\usepackage{color}
\usepackage{colortbl}
\usepackage{rotating}
\usepackage{multirow}
\usepackage{longtable}
\usepackage{capt-of}
\usepackage[font=small]{caption}
\usepackage[round]{natbib}
\usepackage{apalike}
\usepackage{url}
\usepackage{rotating}
\usepackage{soulutf8}

\newcolumntype{Y}{>{\centering\arraybackslash}X}

\definecolor{maroon}{cmyk}{0,0.87,0.68,0.32}
\definecolor{lightblue}{rgb}{.68,.85,.9}

\newcommand{\dem}{\par \noindent{\bf Proof:} \par}
\newcommand{\fin}{\hfill $\square$  \par \bigskip}

\newtheorem{propo}{\bf Proposition}[section]

\newtheorem{rem}{\bf Remark}[section]

\allowdisplaybreaks[1]

\newcommand\LU[1]{\textcolor{blue}}

\usepackage{booktabs} % To thicken table lines

\newcommand{\pcp}{\ensuremath{p}CP\xspace}
\newcommand{\spcp}{S\ensuremath{p}CP\xspace}

\allowdisplaybreaks[1]

\begin{document}
\title{The stratified $p$-center problem}
\author{ Maria Albareda-Sambola$^{1}$, Luisa I. Mart\'{\i}nez-Merino$^2$, Antonio M.\ Rodr\'{\i}guez-Ch\'{\i}a$^2$}

\date{\scriptsize
$^1$ Departamento de Estad\'{\i}stica e
	Investigaci\'on Operativa, Universitat Polit\`{e}cnica de Catalunya-BarcelonaTech, Barcelona, Spain \\
$^2$Departamento de Estad\'{\i}stica e Investigaci\'on Operativa, Universidad de C\'adiz, Spain \today}

\maketitle
\begin{abstract}
	This work presents an extension of the $p$-center problem. In this new model, called Stratified $p$-Center Problem  (S$p$CP),  the demand is concentrated in a set of sites and the population of these sites is divided  into different strata depending on the kind of  service that they require. The aim is to locate $p$ centers to cover the different types of services demanded minimizing the weighted average of the largest distances associated with each of the different strata. In addition, it is considered that more than one stratum can be present at each site. Different formulations, valid inequalities and preprocessings are developed and compared for this problem. An application of this model is presented in order to implement a heuristic approach based on the Sample Average Approximation method (SAA) for solving the probabilistic $p$-center problem in an efficient way.
\end{abstract}


\section{Introduction}
Discrete location problems have been widely studied since the seminal paper  \st{of} \cite{MIP1}, where the first MILP formulation for such a problem was proposed.  Among the fundamental problems in this area, the $p$-Center Problem (\pcp) aims at selecting, from $n$ given sites, the locations of $p$ service centers that minimize the maximum distance between any of the sites and its closest service center. This model, in contraposition to the $p$-median problem, was motivated by the need not to discriminate spatially dispersed clients when locating essential or emergency centers \citep[see][for more details]{garfinkelPC1,calik2015p}.  In the last decades, several extensions of the \pcp have been studied in the literature. These include variants considering capacities \citep{capacPCP,capacPCP2}, pre-existing centers, as in the conditional $p$-center problem \citep{conditional}, or problems with uncertain parameters, which have been addressed both, from the perspective of robust optimization \citep{robust2,robust} and of stochastic programming \citep{SPpcenter,StocCap,PpCP}.

A common characteristic of most of the considered problem variants is that customers are assumed to be homogeneous in the sense that they are all considered in the same way in the objective function. The only exception would be the weighted \pcp, where the distances between each site and its closest center are affected by site-dependent weights.

 In this paper, we consider situations where, for instance, the population of a  region is divided into different strata, and people of some of the strata live together in each of the cities. The problem is to locate centers to cover the essential services of these cities. Due to social or political reasons, the evaluation of the service is measured separately for each stratum. This problem is called the Stratified $p$-Center Problem (S$p$CP). Humanitarian relief planning also fits this model, where different needs of the population must be covered from the located centers, and demand for these different needs may be distributed in a spatially different way. This idea has been used in covering problems \citep{schilling1979} but, up to the best of our knowledge, it has not been applied in the context of the \pcp.

The paper is organized as follows. In Section \ref{classic}, a formulation for the S$p$CP based on the \cite{daskin1995} and \cite{Calik} formulations for the \pcp is introduced. In Section \ref{covering}, some alternative formulations together with some valid inequalities are proposed. Section \ref{saa} applies the results in previous sections for an efficient implementation of a Sample Average Approximation heuristic for the probabilistic $p$-center problem \citep[see][]{PpCP}.
 Section \ref{computational} reports the computational results comparing all the proposed formulations and the results of using Sample Average Approximation. Finally, Section \ref{conclusion} gathers the conclusions of the paper.

\section{Notation and classical formulation}\label{classic}

Let $N=\{1,\ldots,n\}$ be a given set of sites and $p\geqslant 2$ the number of facilities to be located. For each pair $i,j\in N$, let $d_{ij}$ be
the distance from location $i$ to $j$. Besides, $d_{ii}=0$ for $i\in N$ and $d_{ij}>0$ for $i\neq j$.
In the following we use the next notation. The sorted distances associated with pairs of sites are denoted by $$0=d_{(1)}<d_{(2)}<\ldots<d_{(G)}.$$
The sorted distances from a site $i\in N$ to the remaining sites are denoted  by $$0=d_{i(1)}<d_{i(2)}<\ldots<d_{i(G_i)}.$$
In the previous notation $G$ and $G_i$ are the number of different distances between pairs of sites and between $i$ and any other site, respectively (removing possible multiplicities).

Moreover, the population of each site $1,\ldots,n$ is partitioned into a set of strata, taking into account that more than one stratum can be present in a site $i$ and not always all the strata are present in a site.  Given ${\cal S}$ the set of strata in which the population is divided, we consider a family of subsets $\{N^s\}_{s\in \mathcal{S}}$
such that $N^s\subseteq N$ is the set of sites where stratum $s$ is present for $s\in\mathcal{S}$. Then, the sorted distances from a stratum, i.e., the sorted sequence of family $\{d_{ij}\}_{i\in N^s,j\in N}$ is denoted by
$$0=d_{(1)}^s<d_{(2)}^s<\ldots<d_{(G^s)}^s,$$
\noindent where $G^s$ is the number of different distances
of the family $\{d_{ij}\}_{i\in N^s,j\in N}$.

 The problem addressed in this work is based on the classical \pcp. However, in contrast with the \pcp, this new problem considers that population of the sites is divided in different strata depending on the kind of  service that they require. For a given stratum $s$, this problem takes into account the largest distance from the sites where stratum $s$ is present and their corresponding closest service facility. Recall that in the same site there can be inhabitants belonging to more than one stratum.

For each site $j\in N$, and each stratum $s\in\mathcal{S}$, the following binary parameter is defined:
$$ \xi_j^s=
\begin{cases}
1,& \mbox{ if $j\in N^s$,}\\
0,&\mbox{otherwise.}
\end{cases}
$$

Besides, each stratum has an associated weight, ($w_s,\,s\in \mathcal{S}$) that is used to balance the cost related to the different strata in the objective function. The weights can be interpreted in different ways. For instance, they can measure the importance given to a certain stratum.

Given the former parameters, the aim of this problem is to locate $p$ service facilities minimizing the weighted sum of the largest assignments within each stratum. Therefore, the problem can be expressed synthetically in the following way:
 \begin{equation}
\label{probdef}
\min_{\stackrel{P\subseteq N}{|P|=p}}\sum_{s\in\mathcal{S}}w_sd(P, N^s),
\end{equation}
\noindent where  $\displaystyle d(P, N^s)=\max_{j\in  N^s}\min_{i\in P}d_{ij}$. For a given site $j\in N$, we will refer to $\displaystyle\min_{i\in P}d_{ij}$ as the allocation distance of site $j$, so $d(P,N^s)$ is the maximum allocation distance among the sites with presence of stratum $s$, or equivalently within stratum $s$.

The problem previously described can be formulated using the classic $p$-center formulation \citep[see][]{daskin1995}. With this purpose, the following variables are defined:
\begin{eqnarray}
	x_{ij}&=&
	\begin{cases}
		1,&\mbox{if site $j$ is assigned to center $i$,}\\
		0,&\mbox{otherwise,}
	\end{cases}\mbox{ for $i,j\in N$.}\\
	\theta^s&=&\mbox{largest allocation distance for the sites where stratum $s$ is present, $s\in \mathcal{S}$.\label{thetavar}}
\end{eqnarray}
Using these variables, the derived formulation is,
\begin{alignat}{3}
\mbox{(F1)}\quad\quad\quad &&\min\ & \sum\limits_{s\in \mathcal{S}} w_s\theta^s &&\\
&&\textrm{s.t. } & \sum\limits_{i\in N} x_{ii} = p, &&\label{21}\\
&&& \sum\limits_{i\in N} x_{ij} =1, &\qquad& j\in N,\label{22}\\
&&& x_{ij} \leqslant x_{ii}, && i,j\in N,\label{23}\\
&&& \theta^s \geqslant \sum\limits_{i\in N} d_{ij}x_{ij}, && s\in \mathcal{S}, j\in N^s
,\label{24}\\
&&& x_{ij} \in \{0, 1\}, &&i,j\in N,\\
&&& \theta^s \geqslant 0, && s\in \mathcal{S}.
\end{alignat}

Constraint \eqref{21} restricts that there are $p$ centers. Constraints \eqref{22} indicate that each site is associated with only one center. Constraints \eqref{23} restrict that sites must be assigned to an open center. Constraints \eqref{24} ensure that the largest allocation distance within stratum $s$ is not smaller than the allocation distance of any site where stratum $s$ is present. As mentioned before, the objective function is the weighted sum of the largest distances within each stratum.

To the best of our knowledge, the most recent formulation for the $p$-center problem was given by \cite{Calik} providing very good results. We propose a formulation of our problem inspired in \cite{Calik} using the following families of variables.
\begin{eqnarray*}
	\bar{u}_{sr}&=&
	\begin{cases}
		1,&\mbox{ if $d_{(r)}$ is the largest allocation distance among the sites in $N^s$,}\\
		0,&\mbox{otherwise},
	\end{cases}\\
&&\mbox{ $s\in \mathcal{S},r=1,\ldots,G$.}\\
y_i&=&
\begin{cases}
	1,&\mbox{if a center is placed at $i$,}\\
	0,&\mbox{otherwise,}
\end{cases}\mbox{ for $i\in N$.}
\end{eqnarray*}

Using these families of variables, the new formulation is given by

\begin{alignat}{3}
\mbox{(F2)}\quad\quad\quad&&\min\ & \sum\limits_{s\in \mathcal{S}}\sum\limits_{k=1}^G w_sd_{(k)} \bar{u}_{sk} \label{obj1}&&\\
&&\textrm{s.t. } & \sum\limits_{i\in N}y_i = p, &&\label{11}\\
&&& \sum\limits_{k=1}^G \bar{u}_{sk}=1, &&s\in \mathcal{S},\label{12}\\
&&&\sum_{k'=1}^{k-1} \bar{u}_{sk'}\leqslant \sum\limits_{\stackrel{i\in N}{d_{ij}< d_{(k)}}}y_i,&\qquad&  %(j,s)\in N\times\mathcal{S}: \xi_j^s=1
s\in \mathcal{S}, j\in N^s,k=2,\ldots,G,\label{F3}\\
&&& y_i \in \{0, 1\}, &&i\in N,\\
&&& \bar{u}_{sk}\in\{0,1\}, &&s\in \mathcal{S},k=1,\ldots,G.
\end{alignat}

Constraint \eqref{11} restricts that there are $p$ centers. Constraints \eqref{12} ensure that for each stratum, only one of the distances is the largest allocation distance. Constraints \eqref{F3} determine that the largest allocation distance within a stratum $s$ will be among the first $k$ distances if there is a center with a distance smaller than or equal to $d_{(k)}$ with respect to any site in $N^s$.

Observe that $\bar{u}$-variables determine the largest allocation distance
among the sites where each stratum $s\in\mathcal{S}$ is present. As a consequence, only the distances associated with sites in $N^s$
will be necessary to obtain the largest distance with respect to $s$. Therefore, the number of variables can be reduced defining $\tilde{u}$-variables in the following way,
$$
	\tilde{u}_{sk}=
	\begin{cases}
		1,&\mbox{ if $d_{(k)}^s$ is the largest allocation distance for the sites in $N^s$}\\
		0,&\mbox{otherwise},
	\end{cases} s\in \mathcal{S},k=1,\ldots,G^s.
$$

Observe that in the original formulation F2, the number of $\bar{u}$-variables is $|\mathcal{S}| G$. However, by doing this reduction, the obtained number of variables is $\displaystyle\sum_{s\in \mathcal{S}} G^s$.
Taking advantage of this reduction of the number of variables, the new objective function for the model is
\begin{equation}
\sum\limits_{s\in \mathcal{S}}\sum\limits_{k=1}^{G^s} w_sd_{(k)}^s \tilde{u}_{sk},\label{newobj}
\end{equation}
and constraints \eqref{F3} can be replaced by
\begin{equation}
\sum_{k'=1}^{k-1} \tilde{u}_{sk'}\leqslant \sum\limits_{\stackrel{i\in N}{d_{ij}< d_{(k)}^s}}y_i, \quad s\in \mathcal{S}, j\in N^s
, k=2,\ldots,G^s.\label{F31_3}
\end{equation}
 Therefore, this new family of $\tilde{u}$-variables allows us to provide a new formulation with a smaller number of variables and constraints. Moreover, the following result allows to strengthen this new formulation.
\begin{propo}\label{propo21} For $s\in \mathcal{S}$ and $j\in N^s$, let $l_{jr}^s\in\{1,\ldots,G^s\}$ be such that $d_{j(r)}=d_{(l_{jr}^s)}^s$. Considering formulation F2 with $\tilde{u}$ variables (instead of $\bar{u}$ variables), the objective function \eqref{newobj} and replacing \eqref{F3} by
	\begin{equation}
	\sum_{k'=1}^{l_{jr}^s-1} \tilde{u}_{sk'}\leqslant \sum\limits_{\stackrel{i\in N}{d_{ij}< d_{(l_{jr}^s)}^s}}y_i, \quad s\in \mathcal{S}, j\in N^s,r= 2,\ldots,G_j,\label{F31_3n}
	\end{equation}
	results in a valid equivalent formulation F2' with a smaller number of constraints.
	\label{propF3}
\end{propo}
\dem
We prove that constraint families \eqref{F31_3} and \eqref{F31_3n} are equivalent.
Let $\tilde{s}\in \mathcal{S}, \tilde{\j}\in N^s$
 and $\tilde{r}\in\{2,\ldots,G_{\tilde{\j}}\}$.
 Consider the  following subset of constraints of family \eqref{F31_3},
\begin{equation}
\sum_{k'=1}^{ k-1}\tilde{u}_{\tilde{s}k'}\leqslant\sum_{\stackrel{i\in N}{d_{i\tilde{\j}}< d_{(k)}^{\tilde{s}}}}y_i,\quad\quad\quad k\in\{l_{\tilde{\j},\tilde{r}-1}^{\tilde{s}}+1,\ldots,l_{\tilde{\j}\tilde{r}}^{\tilde{s}}\}.\label{subset1}
\end{equation}
Observe that $\displaystyle\sum_{\stackrel{i\in N}{d_{i\tilde{\j}}< d_{(l_{\tilde{\j},\tilde{r}-1}^{\tilde{s}}+1)}^{\tilde{s}}}}y_i=\ldots=\sum_{\stackrel{i\in N}{d_{i\tilde{\j}}< d_{(l_{\tilde{\j}\tilde{r}}^{\tilde{s}})}^{\tilde{s}}}}y_i$, then since
\begin{equation*}
\sum_{k'=1}^{l_{\tilde{\j}\tilde{r}-1}^{\tilde{s}}}\tilde{u}_{\tilde{s}k'}\leqslant\ldots\leqslant\sum_{k'=1}^{l_{\tilde{\j}\tilde{r}}^{\tilde{s}}-1}\tilde{u}_{\tilde{s}k'},
\end{equation*}
the family of constraints \eqref{subset1} is dominated by
$\displaystyle\sum_{k'=1}^{l_{\tilde{\j}\tilde{r}}^{\tilde{s}}-1}\tilde{u}_{\tilde{s}k'}\leqslant \sum_{\stackrel{i\in N}{d_{i\tilde{\j}}<d_{(l_{\tilde{\j}\tilde{r}}^{\tilde{s}})}^{\tilde{s}}}}y_i$.

Therefore, the obtained formulation F2' is equivalent to F2 with less  constraints. In fact, the number of constraints \eqref{F31_3} is $\displaystyle\sum_{j\in N}\sum_{s\in\mathcal{S}}\xi_j^sG^s$ 
and the number of constraints \eqref{F31_3n} is $\displaystyle\sum_{j\in N}\sum_{s\in\mathcal{S}}\xi_j^sG_j$. It is straightforward that for each pair, $s\in \mathcal{S}, j\in N^s$
, $G_j\leqslant G^s$ since, at least, the distances associated with location $j$ must be among the distances related to stratum $s$.

\fin

\section{Formulation using covering variables}\label{covering}
\subsection{ Formulation with stratum-covering variables}

In this subsection we present a formulation based on the use of $y$-variables described in the previous section and the following family of variables:
\begin{eqnarray*}
	u_{sk}&=&
	\begin{cases}
		1,& \mbox{if the largest allocation distance for the sites in $N^s$
is at least $d_{(k)}^s$},\\
		0,&\mbox{otherwise,}
	\end{cases}\\
	&&\mbox{for $s\in \mathcal{S}$, $k=2,\ldots,G^s$.}
\end{eqnarray*}
 Observe that we have used the same strategy as in the former section, so that for each $s\in \mathcal{S}$ the number of $u$ variables will be equal to the number of different distances associated with $s$.
The use of this type of variables for the classical \pcp was introduced by \cite{Elloumi}. Inspired in this idea, we provide the following formulation for the \spcp.
\begin{alignat}{3}
\mbox{(F3)}\quad\quad\quad &&\min\ &\sum_{s\in \mathcal{S}}w_s\left(\sum_{k=2}^{G^s}(d_{(k)}^s-d_{(k-1)}^s)u_{sk}\right) \label{obj}&&\\
&&\textrm{s.t. } & \sum\limits_{i\in N} y_{i} =p, &&\label{31}\\
&&&
u_{sk} \geqslant 1-\sum_{\stackrel{i\in N}{d_{ij}< d_{(k)}^s}}y_i,&& s\in \mathcal{S}, j\in N^s,\, k=2,\ldots,G^s,\label{33}\\
&&& y_{i} \in \{0, 1\}, &&i\in N,\label{aux1}\\
&&& u_{sk}\in \{0,1\},&&s\in\mathcal{S}, k=2,\ldots,G^s.\label{aux2}
\end{alignat}
As it can be seen in \eqref{obj}, the objective function for this formulation can be expressed using a telescopic sum. Constraint \eqref{31} ensures that there are $p$ open centers. Constraints \eqref{33} determine that if there is not a center at a distance smaller than $d_{(k)}^s$ from a site $j\in N^s$, then $u_{sk}=1$.

\begin{propo}\label{propo31}
	Replacing \eqref{33} in F3 by the following families of constraints
	\begin{alignat}{2}
u_{s,l_{jr}^s}&\geqslant 1-\sum_{\stackrel{i\in N}{d_{ij}< d_{j(r)}}}y_i,& \qquad&
s\in \mathcal{S}, j\in N^s,\,r=2,\ldots,G_j,\label{F22}\\
	u_{sk}&\leqslant  u_{s,k-1},  &&s\in\mathcal{S},\, k=3,\ldots,G^{s},\label{F22_2}
	\end{alignat}
	results in an equivalent formulation, F3-\eqref{33}+\eqref{F22}+\eqref{F22_2}.
\end{propo}

\dem
Let $(\tilde{\j},\tilde{s})\in N\times\mathcal{S}$ such that $\xi_{\tilde{\j}}^{\tilde{s}}=1$ and $\tilde{r},\tilde{r}+1\in \{2,\ldots, G_{\tilde{\j}}\}$. Consider the following subset of constraints of family \eqref{33},
\begin{equation}
u_{\tilde{s}k}\geqslant 1-\sum_{\stackrel{i\in N}{d_{i\tilde{\j}}<d_{(k)}^{\tilde{s}}}}y_i,\,\,\,k\in\{l_{\tilde{\j},\tilde{r}-1}^{\tilde{s}}+1,\ldots,l_{\tilde{\j}\tilde{r}}^{\tilde{s}}\}.\label{subset2}
\end{equation}

Observe that {$\displaystyle \sum_{\stackrel{i\in N}{d_{i\tilde{\j}}<d_{(l_{\tilde{\j},\tilde{r}-1}^{\tilde{s}}+1)}^{\tilde{s}}}}y_i=\ldots= \sum_{\stackrel{i\in N}{d_{i\tilde{\j}}<d_{(l_{\tilde{\j}\tilde{r}}^{\tilde{s}})}^{\tilde{s}}}}y_i$, then using \eqref{F22_2}, the family \eqref{subset2} is dominated by
\begin{equation*}
u_{s,l_{\tilde{\j}\tilde{r}}^{\tilde{s}}}\geqslant 1-\sum_{\stackrel{i\in N}{d_{i\tilde{\j}}<d_{(l_{\tilde{\j}\tilde{r}}^{\tilde{s}})}^{\tilde{s}}}}y_i.
\end{equation*}\fin
\begin{rem}
	Formulation F3-\eqref{33}+\eqref{F22}+\eqref{F22_2} has a smaller number of constraints than F3 if $$\sum_{s\in \mathcal{S}}\left(\sum_{j\in N}\xi_{j}^s (G^s-G_j)-G^s+2\right)\geqslant 0. $$
\end{rem}

\subsection{ Formulation with site-covering variables} \label{zvar}

In this section we propose a new formulation for our problem using the following set of variables, inspired in the ones defined by \cite{Garcia} for the \pcp:
\begin{eqnarray*}
	z_{ir}&=&\begin{cases}
		1,& \mbox{if the allocation distance of site $i$ is at least $d_{i(r)}$,}\\
		0,&\mbox{otherwise,}
	\end{cases}\mbox{ for $i\in N$, $r=2,\ldots,G_i$.}
\end{eqnarray*}
Based in this set of variables and $\theta^s$-variables defined by \eqref{thetavar}, we propose the following formulation for our problem:
\begin{alignat}{3}
\mbox{(F4)}\quad\quad\quad &&\min\ &\sum_{ s\in {\cal S}} w_s\theta^s &&\nonumber\\
%&&\textrm{s.t. } & \sum\limits_{i\in N}(1-z_{i2}) =n- p,&&\label{41}\\
&&\textrm{s.t. } &  \sum\limits_{i\in N}z_{i2} =n-p,&&\label{41}\\
&&& \sum_{\stackrel{i\in N}{d_{ij}< d_{j(r)}}}(1-z_{i2})\geqslant 1-z_{jr},&\quad&j\in N,r=3,\ldots,G_j\label{42}\\
&&&\theta^s\geqslant d_{j(r)}z_{jr},&& s\in \mathcal{S}, j\in N^s,r=2,\ldots,G_j,\label{43}\\
&&&  z_{jr}\in \{0,1\},&&  j\in N,r=2,\ldots,G_j,\\
&&&\theta^s\geqslant 0,&&s\in\mathcal{S}.
\end{alignat}
Constraint \eqref{41} indicates that there are $p$ centers. Constraints \eqref{42} ensure that if $z_{jr}=0$ then, there is at least one center at $i$ with $d_{ij}<d_{j(r)}$, i.e., location $j$ is served by a center at a distance smaller than $d_{j(r)}$. Finally, constraints \eqref{43} ensure that $\theta^s$ is the largest allocation distance for sites in $N^s$.

	\begin{propo}
		Formulation F4 is still valid after relaxing the integrality of variables $z_{ir}$ for $i\in N$, $r=3,\ldots,G_i$.	
	\end{propo}
	\dem
		Let $(\tilde{\theta},\tilde{z})$ be an optimal solution of F4 relaxing $z_{ir}$ for $i\in N$, $r=3,\ldots,G_i$.
  We distinguish between two cases.
	
	If $\displaystyle\sum_{\stackrel{i\in N}{d_{i,i_0}<d_{i_0 r_0}}}(1-\tilde{z}_{i2})=0$ then $\tilde{z}_{i_0 r_0}\geqslant 1$ due to constraints \eqref{42}.  Therefore, $\tilde{z}_{i_0 r_0}=1$.
	
	If $\displaystyle\sum_{\stackrel{i\in N}{d_{i,i_0}<d_{i_0 r_0}}}(1-\tilde{z}_{i2})\geqslant 1$, then constraints \eqref{42} reduce to $z_{i_0r_0}\geqslant 0$. Since positive values of  $\tilde{z}_{i_0r_0}$ penalize the objective function due to constraints \eqref{43}, then $\tilde{z}_{i_0r_0}=0$.
	\fin
	
	Preliminary computational results show that this relaxation does not improve computational times of formulation F4.

\begin{propo}
Replacing constraints \eqref{43} in F4 by
\begin{equation}
\theta^s\geq\sum_{r=2}^{G_j}(d_{j(r)}-d_{j(r-1)})z_{jr},\quad s\in \mathcal{S}, j\in N^s
,\label{nuevaF4}
\end{equation}
results in a valid formulation F4-\eqref{43}+\eqref{nuevaF4} for the problem with less constraints, that dominates F4.
\end{propo}
\dem
Let $s\in \mathcal{S}, j\in N^s$.
Note that, due to constraints \eqref{42} and constraints \eqref{nuevaF4} it holds that $z_{jr}\leqslant z_{j,r-1}$ for $r\in\{3,\ldots,G_j\}$ since,
$$\sum_{\stackrel{i\in N}{d_{ij}<d_{j(r)}}}(1-z_{i2})\geqslant\sum_{\stackrel{i\in N}{d_{ij}<d_{j(r-1)}}}(1-z_{i2}),$$
and $z$-variables penalize in the objective function through constraints \eqref{nuevaF4}. Hence, since $z_{jr}\in\{0,1\}$ we have that
\begin{equation*}
\theta^s=\max_{j\in N^s}\left\{\sum_{r=2}^{G_j}(d_{j(r)}-d_{j(r-1)})z_{jr}\right\}
\end{equation*}
and then the formulation  F4-\eqref{43}+\eqref{nuevaF4} is valid. Moreover, for the relaxed problem we have that
 $$\sum_{r=1}^{G_{j}}(d_{j(r)}-d_{j(r-1)})z_{jr}\geqslant \max_{r=1,\ldots,G_j}d_{j(r)}z_{jr},\,\,\forall s\in \mathcal{S}, j\in N^s,$$
 i.e., this formulation dominates F4. Besides, the number of constraints \eqref{43} is $\displaystyle \sum_{j\in N}\sum_{s\in\mathcal{S}}\xi_j^sG_j$ and the number of constraints \eqref{nuevaF4} is $\displaystyle\sum_{j\in N}\sum_{s\in\mathcal{S}}\xi_j^s$. Then, formulation F4-\eqref{43}+\eqref{nuevaF4} has a smaller number of constraints than F4.
\fin

We have also studied alternative formulations using a non-cumulative version of the $z$-variables, i.e., defining
\begin{eqnarray*}
	\bar{z}_{ir}&=&\begin{cases}
		1,& \mbox{if the allocation distance of site $i$ is $d_{i(r)}$,}\\
		0,&\mbox{otherwise.}
	\end{cases}\mbox{ for $i\in N$, $r=2,\ldots,G_i$.}
\end{eqnarray*}
Nevertheless, a preliminary computational analysis of these formulations shows a worse performance with respect to F4.

\subsection{ Formulation with stratum- and site-covering variables}\label{zuvar}

The last formulation that we propose combines two families of covering variables, one associated with the distances from each stratum $s\in \mathcal{S}$ ($u$-variables) and another one with the allocation of each site $i\in N$ ($z$-variables). The combination of both families of variables is inspired in the formulation of \cite{Marin} for the Discrete Ordered Median problem.

For each $s\in\mathcal{S}$, $k\in\{2,\ldots,G^s\}$ and $i\in N$ we define
$$
\bar{l}_{ik}^s=
\begin{cases}
r,&\mbox{if $r\in\{1,\ldots,G_i\}$ exists such that $d_{i(r)}=d_{(k)}^s$ and $\xi_i^s=1$,}\\
0,&\mbox{otherwise.}\\
\end{cases}
$$

Then, the obtained formulation is
\begin{alignat}{3}
\mbox{(F5)}\quad\quad &&\omit\rlap{$\min\ \sum\limits_{s\in {\cal S}} \sum\limits_{k=2}^{G^{s}}w_s(d_{(k)}^s-d_{(k-1)}^s )u_{sk}$}\\
&&\textrm{s.t. } & \eqref{41},\eqref{42},&&\nonumber\\
&&& u_{sk}\geqslant z_{i,\bar{l}_{ik}^{s}},&\qquad\quad & s\in \mathcal{S}, i\in N^s, k=2,\ldots G^s: \bar{l}_{ik}^s>0,\label{F55}\\
&&&u_{s,k-1}\geqslant u_{sk},&&s\in\mathcal{S},k=3,\ldots,G^{s},\label{F55_3}\\
&&&u_{sk}\in\{0,1\},&&s\in\mathcal{S},k=2,\ldots,G^s,\\
&&&z_{ir}\in\{0,1\},&&i\in N,r=2,\ldots, G_i.
\end{alignat}
 Constraints \eqref{F55} determine the largest allocation distance among the sites in $N^s$. Observe that constraints \eqref{F55_3} are valid inequalities for formulation F5. Indeed, if in a particular solution $u_{sk}=b$ and $u_{s k-1}=a$ with $b>a$, then, a feasible solution with lower objective value can be found by taking $u_{sk}=a$. Constraints \eqref{F55_3} are included in the formulation from the beginning since they provided good results in a preliminary computational study.

Note that constraints \eqref{F55} can be equivalently written in the following way,
\begin{equation}
u_{s,l_{ir}^s}\geqslant z_{ir}\quad\quad s\in \mathcal{S}, i\in N^s,r=2,\ldots,G_i.\label{F2_3}
\end{equation}
Where $l_{ir}^s$ is the index already defined in Proposition \ref{propF3}.
 To derive another valid formulation from (F5), we include the following notation,
$$
l_{ik}^{'s}=\begin{cases}
\min\{r:d_{i(r)}\geqslant d_{(k)}^s\},&\mbox{ if }d_{(k)}^s\leqslant d_{i(G_i)}\\
G_{i}+1,&\mbox{otherwise.}
\end{cases}
$$

\begin{propo}
	By replacing \eqref{F55} in F5 by
	\begin{equation}
	u_{sk}\geqslant z_{i,l_{ik}^{'s}},\quad\quad s\in \mathcal{S}, i\in N^s,k=2,\ldots,G^s,l_{ik}^{'s}\leqslant G_i.\label{desagregada}
	\end{equation}
	a valid formulation, F5-\eqref{F55}+\eqref{desagregada}, with a larger number of constraints is obtained.
\end{propo}
\dem
First, formulation F5-\eqref{F55}+\eqref{desagregada} is valid, since \eqref{desagregada} determine the largest allocation distance among the sites where stratum $s$ is present.

Observe that family of constraints \eqref{F55} is a subset of constraints \eqref{desagregada}  since $l^{'s}_{ik}=\bar{l}_{ik}^s$ when $d_{i(r)}=d_{(k)}^s$ for some $r\in\{2,\ldots,G_i\}$  and $\xi_i^s=1$.
 Therefore F5-\eqref{F55}+\eqref{desagregada} dominates formulation F5. Concretely, the number of constraints \eqref{desagregada} is $\displaystyle\sum_{i\in N}\sum_{s\in\mathcal{S}}\xi_i^s(G^s-1)$. The number of constraints \eqref{F55} is $\displaystyle\sum_{i\in N}\sum_{s\in\mathcal{S}}\xi_i^s(G_i-1)$. As stated before $G_i\leqslant G^s$ for $s\in \mathcal{S}, i\in N^s$. Consequently, the number of constraints \eqref{desagregada} is larger than the number of constraints \eqref{F55}.

\fin

\begin{propo}\label{agregadas}
	\begin{itemize}
		\item[i)] Constraints \eqref{F55} can be replaced by their following aggregated form:
	\begin{equation}
		n_{sk}u_{sk}\geqslant \sum_{\stackrel{i\in N^s}{\bar{l}_{ik}^s\neq 0}}z_{i \bar{l}_{ik}^{s}}\quad s\in\mathcal{S},k=2,\ldots,G^s,\label{F52}
		\end{equation}
		where $n_{sk}=|\{i\in N^s \mbox{ and there exists}$ $r\in\{2,\ldots,G_i\}\mbox{ such that } d_{i(r)}=d_{(k)}^s\}|$. This yields the new valid formulation, F5-\eqref{F55}+\eqref{F52}.
		\item[ii)] Constraints \eqref{desagregada} can be replaced by their aggregated form that can be expressed as
\begin{equation}
n_su_{sk}\geqslant\sum_{\stackrel{i\in N^s}{{l'}_{ik}^s\leqslant G_i}}z_{i,l^{'s}_{ik}},\quad s\in\mathcal{S},k=2,\ldots, G^{s},\label{53}\\
\end{equation}
		Where $n_s=|N^s|$. This yields the new valid formulation F5-\eqref{F55}+\eqref{53}.
	\end{itemize}
\end{propo}
\dem
\begin{itemize}
	\item[i)]
	Observe  that, by \eqref{F52}, variables $u_{sk}$ take the value $1$ if the maximum distance among the sites in $N^s$ is at least $d_{(k)}^s$.
	Indeed, if this allocation distance is at least $d_{(k)}^s$ then, by \eqref{42}, there exists a site $j\in N^s$ such that $z_{j\bar{l}_{jk}^s}=1$ and then, by \eqref{F52}, $u_{sk}=1$.
	
	Moreover, \eqref{F52} are valid since $n_{sk}$ is the maximum value that the right hand side of constraints \eqref{F52} can take.

	\item[ii)] By an argument analogous to the one discussed in  i), we have that formulation F5-\eqref{F55}+\eqref{53} is valid for the \spcp.
\end{itemize}
\fin

Besides, another aggregated version of constraints  \eqref{F55} is:
\begin{equation}
\sum_{\stackrel{s\in \mathcal{S}:\xi_i^s=1}{l_{ir}^s\geqslant 2}}u_{s,l_{ir}^s}\geqslant {\left(\sum_{s\in{\cal S}}\xi_i^s\right)} z_{ir}\quad i\in N,r=2,\ldots,G_i.\label{F6}
\end{equation}
 Some computational studies have been carried out with formulation  F5-\eqref{F55}+\eqref{F6}. However, it provides worse running times that formulations presented in Proposition \ref{agregadas}.

	\begin{propo}\label{relaj_z}
		Formulation F5 and all its variants (F5-\eqref{F55}+\eqref{F2_3}, F5-\eqref{F55}+\eqref{desagregada}, F5-\eqref{F55}+\eqref{F52}, F5-\eqref{F55}+\eqref{53}, F5-\eqref{F55}+\eqref{F6}) remain valid if integrality of variables $z_{ir}$ is relaxed for $i\in N$, $r=3,\ldots,G_i$.
	\end{propo}
\dem
Let $(\tilde{u},\tilde{z})$ be an optimal solution of the model relaxing $z_{ir}$ for $i\in N$, $r=3,\ldots,G_i$. If $(\tilde{u},\tilde{z})$ are all binary, we are done. Otherwise, there is at least one $0 < \tilde{z}_{i_0 r_0} < 1$ with $i_0\in N$, $r_0\in\{3,\ldots,G_i\}$. For this variable, Constraint \eqref{42} reduces to $z_{i_0r_0}\geqslant 0$, since $z_{i2}$ are binary for $i\in N$. Hence $\tilde{z}_{i_0 r_0}$ value can be replaced by $0$ without violating these constraints. Besides, constraints \eqref{F55}, \eqref{F2_3}, \eqref{desagregada}, \eqref{F52}, \eqref{53} or \eqref{F6} (depending on the variant of F5) are not violated if $\tilde{z}_{i_0r_0}$ takes value $0$ and the objective value is not worse.
\fin

Computational results in Section \ref{computational} show that this relaxation improves the times of formulation F5-\eqref{F55}+\eqref{53}.
	
\begin{propo} \label{proprorelax}
		Formulations F5, F5-\eqref{F55}+\eqref{F2_3} and F5-\eqref{F55}+\eqref{desagregada},  remain valid if we relax the integrality condition of $u_{sk}$  variables for $s\in {\cal S}$, $k=2,\ldots,G^s$ and $z_{ir}$ variables for $i\in N$, $r\in\{3,\ldots,G_i\}$.
	\end{propo}
\dem
Since $z$-variables take integer values as observed in Proposition \ref{relaj_z} and since $u_{sk}$ for $s\in {\cal S}$, $k=2,\ldots,G^s$ penalize the objective function, it holds that $u_{sk}$ take integer values due to constraints \eqref{F55} (or equivalently, due to constraints \eqref{F2_3} or \eqref{desagregada}).\
\fin

 Preliminary computational results show that the relaxations introduced in Proposition \ref{proprorelax} do not improve the running times of the corresponding models.

\subsection{ Reducing the number of covering variables} \label{reduction}

Observe that some of $z$-variables described in formulations of subsections \ref{zvar} and \ref{zuvar} could be fixed. Since $p$ centers are located in the S$p$CP, then the distance associated with a client $i$ will not be among the $p-1$ worst possible ones. Then, the following constraints allow to fix some variables.

Let $\tilde{d}_{i(1)}\leqslant \tilde{d}_{i(2)}\leqslant\ldots\leqslant\tilde{d}_{i \tilde{G}_i}$ be the sorted distances of all possible assignments of site $i$  (observe that this sequence of distances can contain repeated values), then
\begin{equation}
z_{ir}=0\,\,\,\forall i\in N, r\in\{2,\ldots,G_i\} \mbox{ such that }d_{i(r)}>\tilde{d}_{i(n-p+1)}.\label{P3}
\end{equation}

Consequently, for each $i\in N$ it is only necessary to define $z_{ir}$ for $r=2,\ldots,G_i$ such that $d_{i(r)}\leqslant\tilde{d}_{i(n-p+1)}$ .

Regarding $u$-variables appearing in formulations F3 and F5, observe that these are binary variables indicating for each stratum $s\in\mathcal{S}$ whether the largest distance associated with stratum $s$ is at least $d_{(k)}^s$  or not, where $k=2,\ldots, G^s$.
The number of  $u$-variables for each stratum $s\in\mathcal{S}$ is $G^s-1$, i.e, the number of different distances from sites in $N^s$ to all candidate locations (excluding distance $0$). In this subsection, we analyze if the number of $u$-variables can be reduced for each stratum.

In fact, the number of $u$-variables could be reduced if tighter bounds on the largest allocation distance associated with each stratum for the S$p$CP were known. The following proposition exploits this argument.

\begin{propo}
	For each stratum $s\in\mathcal{S}$, let $v(p\mbox{CP}_s)$ be the optimal value of a $p$-center problem where the set of candidates centers is $N$ and the set of demand points is $N^s$, from now on, denoted with $p\mbox{CP}_s$. Then, the largest allocation distance associated with $s$  is at least $v(p\mbox{CP}_s)$ in the optimal solution of the S$p$CP.
\end{propo}

\dem
 Observe that the solution of the S$p$CP is feasible for the $p$CP$_s$. Then, given a solution of S$p$CP, its objective value for $p$CP$_s$ will be greater than or equal to $v(pCP_s)$.
\fin

As a result, if a lower bound or the optimal value of $p$CP$_s$ is obtained, then the number of $u$-variables associated  with stratum $s$ can be reduced.  To reduce the number of variables we can follow the next scheme for each $s\in \mathcal{S}$:

\begin{itemize}
  \item Obtain a lower bound on the $p$CP$_s$ or its optimal objective value. This value can be denoted as $LB_{s}$.
  \item Define $u_{sk}$ variables for all $k\in\{h:2\leqslant h\leqslant G^s\mbox{ and } d_{(h)}^s  > LB_{s}\}$.
  \item  For each $s\in{\cal S}$, given that $d_{(k_s)}^s$ is the largest distance associated with stratum $s$ such that $d_{(k_s)}^s\leqslant LB_s$, the considered objective function will be:
		$$ \sum_{s\in  {\cal S}}w_s \sum_{k=k_s+1}^{G^{s}}\left(d_{(k_s)}^s+(d_{(k)}^s-d_{(k-1)}^s )u_{sk}\right)$$
 Observe that this is equivalent to fix $u_{sk}=1$ for $k\leqslant k_s, s\in {\cal S}$.
		
\end{itemize}

Several criteria can be used to obtain an adequate bound $LB_s$ for each stratum. In particular, in the computational experiments of this work we present  two ways for obtaining these bounds. The first one uses the linear relaxation of the $pCP_s$ using the classic formulation of \cite{daskin1995}. The second  one consists in using the binary algorithm proposed in \cite{Calik}.

Observe that the argument described in constraints \eqref{P3} for $z$-variables could be also useful to fix some of the $u$-variables. In particular, the following variables can be fixed:
\begin{equation}\label{varured}
u_{sk}=0,\,\, (s,k)\in {\cal K},
\end{equation}
where ${\cal K}$ is the set of pairs $(s,k)\in {\cal S}\times \{2,\ldots,G^s\}$ such that for every $i\in N^s$, $d_{(k)}^s>\tilde{d}_{i(n-p+1)}$.

Summing up, $u$-variables can be reduced using the scheme described before and constraints \eqref{varured}. In section \ref{computational} we study the percentage of $z$- and $u$- variables fixed by applying the former criteria.

\subsection{Valid inequalities for F5}\label{inequalitiessection}

Some constraints related to closest assignments could be applied for this problem. Some of the constraints appearing in \cite{cac} have been adapted for formulation F5 (the most promising formulation as we will see in Section \ref{computational}). However, the only valid inequality that presents good results is the one described  below:

\begin{equation}
z_{ir} \leqslant z_{j2}\quad i,j\in N,r=2,\ldots G_i:d_{i(r-1)}=d_{ij},\label{R1mod}
\end{equation}
These constraints could be considered as derived from the ones proposed by \cite{Dobson}. As observed, given $i,j\in N$ these constraints restrict the distance associated with $i$ to be smaller than or equal to the distance $d_{ij}$ if a center is located at $j$.

In the following we introduce other valid inequalities that take advantage of the relationship between two different strata.
\begin{equation}
\sum_{k=2}^{G^{s_1}}(d_{(k)}^{s_1}-d_{(k-1)}^{s_1})u_{s_1k}\leqslant\sum_{k=2}^{G^{s_2}}(d_{(k)}^{s_1}-d_{(k-1)}^{s_1})u_{s_2k},\,\, s_1,s_2\in S:N^{s_1}\subseteq N^{s_2}.\label{escenarios}\\
\end{equation}
These constraints state that  the largest allocation distance associated with stratum $s_1$ will be smaller than or equal to the  one associated with stratum $s_2$ if stratum $s_2$ is present in each site of $N^{s_1}$.
Similarly the next constraints follow:
\begin{equation}
u_{s_1k} \leqslant u_{s_2l}, s_1,s_2\in S,k=2,\ldots,G^{s_1},l=2,\ldots,G^{s_2}:N^{s_1}\subseteq N^{s_2}, d_{(k)}^{s_1}=d_{(l)}^{s_2},\label{escenarios2}
\end{equation}
Constraints \eqref{escenarios2} hold since if the largest allocation distance associated with $s_2$ is smaller than $d_{(l)}^{s_2}$ and $N^{s_1}\subseteq N^{s_2}$, then the largest  allocation distance within $s_1$ cannot be greater than or equal to $d_{(l)}^{s_2}=d_{(k)}^{s_1}$. The accumulated version of these valid inequalities is:
\begin{equation}
\sum_{k=2}^{G^{s_1}}u_{s_1k}\leqslant\sum_{k=2}^{G^{s_2}}u_{s_2k}, s_1,s_2\in S:N^{s_1}\subseteq N^{s_2}.\label{escenarios3}
\end{equation}

Other valid inequalities are  those ensuring that $z$ variables are sorted in non-increasing order for each $i\in N$, i.e.,
\begin{equation}
z_{ir}\geqslant z_{i,r+1},\quad i\in N,r=2,\ldots,G_i-1.\label{Restz}
\end{equation}

 All these valid inequalities will be analyzed in Section \ref{computational}.

\section{SAA for the probabilistic $p$-center problem}\label{saa}
Recall from \cite{PpCP} that the Probabilistic $p$CP (P$p$CP) is defined as the variant
of the $p$CP where sites represent potential demand points, and the locations of the $p$ centers have to be decided before the actual subset of sites that need to be served is revealed. In this problem, the goal is to minimize the expected maximum distance between a site with demand and its closest center. Here, expectation is computed with respect to the probability distribution of the binary random vector defining the subset of sites that have demand.

Notice that, in fact, when uncertainty is modeled by means of a set of scenarios, the P$p$CP can
be cast as a S$p$CP. In this case, each stratum would represent the set of sites having demand at a given scenario, and the stratum weight would correspond to the corresponding scenario probability. This suggests exploiting the S$p$CP formulations presented in this paper to solve the P$p$CP using the well-known Sample Average Approximation method (SAA).

SAA is based on using Monte Carlo Sampling in the probability space defined by the random
variables involved in a problem definition \citep[see][]{Homem}. Although this idea was already used before for solving stochastic programming problems \citep[][]{Rubinstein,Robinson}, the term SAA was formally defined in \cite{Kleiwegt}. We next provide a sketch of this methodology; for more details, see \cite{Shapiro}  or \cite{Linderoth}.

Consider the two stage program (P) $z^*=\min_{x\in X} f(x)+\mathcal{Q}(x)$, where the recourse function is defined as $Q(x)=\mathbb{E}_{\xi}[v(x,\xi)]$ and, given a solution $x$ and a realization of the random vector $\xi$, $\xi_0$,
the so-called second stage problem is $v(x,\xi_0)=\min_{y\in Y(x,\xi_0)}q(y;x,\xi_0) $. Note that if $\xi$ is a discrete random vector with a finite support,
$\Omega$, and each scenario $\omega \in \Omega$
has a known probability $p^{\omega}$, then,
by replicating the variables of the second stage problem, (P) can be equivalently expressed as:
\begin{alignat}{3}
\mbox{($P^{\prime}$) }z^*=&&\min &f(x)+
\sum_{\omega\in\Omega}p^{\omega}q(y^{\omega};x,\xi^{\omega}) &&\\
&&\textrm{s.t. } & x\in X,&&\nonumber\\
&& &y^{\omega}\in Y(x,\xi^{\omega}),&&\omega \in \Omega.\nonumber
\end{alignat}
Accordingly, using a random sample $\Omega^{M}\subset \Omega$, with  $M=|\Omega^M|$, $P$ can be approximated as
\begin{alignat}{3}
\mbox{($P^{M}$) }z^M=&&\min\ &f(x)+\frac{1}{M}\sum_{\omega\in\Omega^M} q(y^{\omega};x,\xi^{\omega}) &&\\
&&\textrm{s.t. } & x\in X,&&\nonumber\\
&& &y^{\omega}\in Y(x,\xi^{\omega}),&&\omega \in \Omega^M.\nonumber
\end{alignat}

Problem $P^M$ is often referred to as sample average approximation problem. It is well known
that given $M$, the expected value of this problem, $\mathbb{E}(z^M)$, is a lower bound on $z^*$ and it converges to $z^*$ as $N$ increases. Moreover, under some mild conditions on $X$ and $v$, the random vector $x^{M,*}$ representing the optimal solution to $P^M$ becomes arbitrarily close to the set of optimal solutions to $P$ with probability 1. A common way to estimate $\mathbb{E}(z^M)$ is to solve a sequence of realizations of $P^M$
for a given sample size $M$, and use the average of the corresponding optimal values as an estimate of $\mathbb{E}(z^M)$. The sequence is evaluated iteratively, and the termination criterion is most often related with the convergence of this average. The best of the solutions obtained in that sequence of problems is kept as a good approximation of the optimal solution.
 The next scheme describes the SAA for case of the P$p$CP.
\begin{itemize}
	\item In each iteration (denoted by $k$):
	\begin{itemize}
		\item {\bf Generation of a random sample $\Omega^k\in\Omega$.} To this aim, a total of $m$ scenarios are generated and each scenario $w=(\xi_1^{\omega},\ldots,\xi_n^{\omega})\in\Omega^k$ is obtained as follows. First for each $i\in N$, create a random number $r\in[0,1)$. If $r<q_i$, where $q_i$ is the probability of client $i$ to have demand, then $\xi_i^{\omega}=1$. Otherwise, $\xi_i^{\omega}=0$. In this case $\xi_i^{\omega}$ determines whether client $i$ has demand in scenario $\omega$ or not.
		\item {\bf Solving of the sample average approximation problem.} Solve the S$p$CP with one of the formulations described in sections \ref{classic} or \ref{covering}. Note that in this case we set ${\cal S}=\Omega^k$.
		\item {\bf Evaluation of the solution.}  Evaluate the S$p$CP solution according to the objective function of the P$p$CP. This allows to obtain an upper bound that can be updated in each iteration if it is improved. Besides, save the optimal objective value of S$p$CP to obtain the average objective value after all the iterations.
	\end{itemize}
	
	\item Stopping  criterion: Stop the procedure after a number of iteration or when the average of the objective values of  \spcp converges. 
\end{itemize}

In the next section we will show some computational results of SAA using random samples of size $M=10$. Besides, we will see how the use of different formulations of S$p$CP can affect the performance of the SAA.

\section{Computational results}\label{computational}

This section is devoted to the computational studies of the formulations described along the paper for the \spcp. The instances used in this computational experience are
based on the $p$-median instances from the ORLIB\footnote{Electronically available at \url{http://people.brunel.ac.uk/~mastjjb/jeb/orlib/files/}}.

 For the smallest instances ($n=6,\ldots,75$), the used matrices are submatrices of instances pmed1, pmed2, pmed3, pmed4 and pmed5 from the ORLIB data. For instances with $n=100,200,300,400$, the matrices are those corresponding to instances pmed1-pmed20. In all cases, several $p$ values are considered ranging between $p=2$ (for the smallest instances) to $p=60$ (for the largest instances). Finally, in Table \ref{tab7} all the ORLIB distance matrices together with their corresponding $p$ values are studied.

For each instance, a total of $|{\cal S}|=10$ strata are generated. Besides, each stratum ($s$) is independently created. First, a number $q_i\in(0,1)$ is associated with each $i\in N$. Then a random number in $r\in[0,1)$ is created. If $r<q_i$, then $\xi_i^{s}=1$. Otherwise, $\xi_i^{s}=0$.

The formulations  are implemented  in the commercial solver Xpress 8.0 using the modeling language Mosel. All the runs  are carried out on the same computer with an Intel(R) Core(TM) i7-4790K processor with 32 GB RAM. We remark that the cut generation of Xpress is disabled to compare the relative performance of formulations cleanly.

First, we report a comparison of all proposed formulations in sections \ref{classic} and \ref{covering}. In this study, we observe that the best results are provided by  a variant of formulation F5. After that, we analyze if valid inequalities and the reduction of variables improve the computational times. Finally, Sample Average Approximation for P$p$CP is implemented using some of the S$p$CP formulations presented before.

\subsection{Comparison of formulations}

Table \ref{tab1} reports the results of the S$p$CP formulations proposed in sections \ref{classic} and \ref{covering}. As can be observed, some formulations include several variants replacing some of the constraints by others. With these new constraints, the aim is to improve the  running times of some of  these formulations. Table \ref{tab1} reports two columns for each formulation. The first one shows the average running time for solving the model and the second column reports the average LP gap (in percentage, $\%$). Note that each entry corresponds to the average over five instances of the same size
and that the reported average running time is the average among the instances that are solved in less than two hours. The number of unsolved instances after two hours is reported in parentheses. Besides, observe that formulation F5-\eqref{F55}+\eqref{53}* corresponds to formulation F5 replacing constraints \eqref{F55} by constraints \eqref{53} relaxing variables $z_{ir}$ for $i\in N$, $r\in\{3,\ldots,G_i\}$.

\begin{sidewaystable}[hp]
	\centering
	\caption{Formulations times and LP gap comparison}
	\scalebox{0.6}{
	\begin{tabular}{rr|lr|lr|rr|lr|rr|rr|rr|rr|rr|rr|rr}
		\hline
			&       & \multicolumn{2}{c|}{F1 } & \multicolumn{2}{c|}{F2} & \multicolumn{2}{c|}{F2'} & \multicolumn{2}{c|}{F3} & \multicolumn{2}{c|}{F3-\eqref{33}+\eqref{F22}+\eqref{F22_2}} & \multicolumn{2}{c|}{ F4} & \multicolumn{2}{c|}{F5} & \multicolumn{2}{c|}{F5-\eqref{F55}+\eqref{desagregada}} & \multicolumn{2}{c|}{F5-\eqref{F55}+\eqref{F52}} & \multicolumn{2}{c}{F5-\eqref{F55}+\eqref{53}}& \multicolumn{2}{|c}{ F5-\eqref{F55}+\eqref{53}*} \\
		\multicolumn{1}{l}{n} & \multicolumn{1}{l|}{p} & Time  & \multicolumn{1}{r|}{LP Gap} & Time  & \multicolumn{1}{r|}{LP Gap} & \multicolumn{1}{r}{Time} & \multicolumn{1}{r|}{LP Gap} & Time  & \multicolumn{1}{r|}{LP Gap} & \multicolumn{1}{r}{Time } & \multicolumn{1}{r|}{LP Gap} & \multicolumn{1}{r}{Time } & \multicolumn{1}{r|}{LP Gap} & \multicolumn{1}{r}{Time } & \multicolumn{1}{r|}{LP Gap} & \multicolumn{1}{r}{Time} & \multicolumn{1}{r|}{LP Gap} & \multicolumn{1}{r}{Time} & \multicolumn{1}{r|}{LP Gap} & \multicolumn{1}{r}{Time } & \multicolumn{1}{r|}{LP Gap} & \multicolumn{1}{r}{Time } & \multicolumn{1}{r}{LP Gap} \\
		\hline
		6     & 2     & \multicolumn{1}{r}{\textbf{0.00}} & 20.96 & \multicolumn{1}{r}{0.02} & 11.12 & 0.01  & 11.12 & \multicolumn{1}{r}{0.02} & 11.12 & 0.02  & 11.12 & 0.01  & 33.77 & 0.02  & 11.12 & 0.02  & 11.12 & 0.01  & 24.20 & 0.01  & 50.73 & 0.01  & 50.73 \\
	10    & 3     & \multicolumn{1}{r}{\textbf{0.02}} & 28.50 & \multicolumn{1}{r}{0.16} & 19.20 & 0.05  & 19.20 & \multicolumn{1}{r}{0.19} & 19.20 & 0.08  & 19.20 & 0.03  & 40.80 & 0.08  & 19.20 & 0.16  & 19.20 & 0.09  & 26.82 & 0.05  & 58.92 & 0.06  & 58.92 \\
	10    & 5     & \multicolumn{1}{r}{\textbf{0.01}} & 44.34 & \multicolumn{1}{r}{0.13} & 28.72 & 0.04  & 28.72 & \multicolumn{1}{r}{0.15} & 28.72 & 0.05  & 28.72 & 0.03  & 52.08 & 0.05  & 28.72 & 0.09  & 28.72 & 0.05  & 43.40 & 0.05  & 74.60 & 0.05  & 74.60 \\
	13    & 3     & \multicolumn{1}{r}{\textbf{0.06}} & 30.20 & \multicolumn{1}{r}{0.89} & 17.60 & 0.18  & 17.60 & \multicolumn{1}{r}{0.61} & 17.60 & 0.23  & 17.60 & 0.09  & 47.82 & 0.23  & 17.60 & 0.56  & 17.60 & 0.26  & 29.86 & 0.12  & 58.12 & 0.11  & 58.12 \\
	13    & 5     & \multicolumn{1}{r}{\textbf{0.04}} & 40.94 & \multicolumn{1}{r}{0.55} & 22.38 & 0.10  & 22.38 & \multicolumn{1}{r}{0.37} & 22.38 & 0.12  & 22.38 & 0.07  & 49.22 & 0.12  & 22.38 & 0.21  & 22.38 & 0.13  & 35.13 & 0.09  & 65.30 & 0.09  & 65.30 \\
	13    & 8     & \multicolumn{1}{r}{\textbf{0.03}} & 48.74 & \multicolumn{1}{r}{0.28} & 24.11 & 0.05  & 24.11 & \multicolumn{1}{r}{0.19} & 24.11 & 0.07  & 24.11 & 0.05  & 51.34 & 0.06  & 24.11 & 0.11  & 24.11 & 0.06  & 39.06 & 0.05  & 78.82 & 0.05  & 78.82 \\
	15    & 3     & \multicolumn{1}{r}{\textbf{0.10}} & 30.14 & \multicolumn{1}{r}{1.39} & 16.08 & 0.31  & 16.08 & \multicolumn{1}{r}{1.40} & 16.08 & 0.47  & 16.08 & 0.14  & 48.61 & 0.38  & 16.08 & 0.75  & 16.08 & 0.40  & 27.60 & 0.17  & 56.42 & 0.18  & 56.42 \\
	15    & 7     & \multicolumn{1}{r}{\textbf{0.06}} & 53.11 & \multicolumn{1}{r}{0.93} & 32.47 & 0.15  & 32.47 & \multicolumn{1}{r}{0.68} & 32.47 & 0.21  & 32.47 & 0.09  & 58.53 & 0.16  & 32.47 & 0.35  & 32.47 & 0.18  & 42.18 & 0.14  & 75.30 & 0.14  & 75.30 \\
	15    & 10    & \multicolumn{1}{r}{\textbf{0.04}} & 46.31 & \multicolumn{1}{r}{0.47} & 19.74 & 0.06  & 19.74 & \multicolumn{1}{r}{0.29} & 19.74 & 0.08  & 19.74 & 0.07  & 49.00 & 0.06  & 19.74 & 0.14  & 19.74 & 0.05  & 38.90 & 0.07  & 83.32 & 0.06  & 83.32 \\
	20    & 3     & \multicolumn{1}{r}{\textbf{0.25}} & 33.59 & \multicolumn{1}{r}{6.92} & 18.87 & 0.94  & 18.87 & \multicolumn{1}{r}{5.12} & 18.87 & 1.27  & 18.87 & 0.31  & 51.65 & 0.85  & 18.87 & 2.43  & 18.87 & 1.32  & 28.76 & 0.47  & 54.20 & 0.51  & 54.20 \\
	20    & 7     & \multicolumn{1}{r}{\textbf{0.21}} & 44.51 & \multicolumn{1}{r}{3.35} & 21.66 & 0.43  & 21.66 & \multicolumn{1}{r}{1.80} & 21.66 & 0.56  & 21.66 & 0.25  & 52.90 & 0.38  & 21.66 & 0.94  & 21.66 & 0.56  & 32.66 & 0.29  & 68.11 & 0.31  & 68.11 \\
	20    & 10    & \multicolumn{1}{r}{\textbf{0.14}} & 52.88 & \multicolumn{1}{r}{2.43} & 27.43 & 0.32  & 27.43 & \multicolumn{1}{r}{1.88} & 27.43 & 0.43  & 27.43 & 0.20  & 58.24 & 0.28  & 27.43 & 0.56  & 27.43 & 0.37  & 40.04 & 0.29  & 77.68 & 0.29  & 77.68 \\
	25    & 3     & \multicolumn{1}{r}{\textbf{0.54}} & 29.89 & \multicolumn{1}{r}{17.66} & 15.68 & 1.88  & 15.68 & \multicolumn{1}{r}{11.86} & 15.68 & 2.19  & 15.68 & 0.57  & 50.16 & 1.23  & 15.68 & 4.06  & 15.68 & 2.34  & 26.84 & 1.02  & 50.37 & 1.13  & 50.37 \\
	25    & 7     & \multicolumn{1}{r}{0.56} & 43.51 & \multicolumn{1}{r}{15.33} & 22.04 & 1.36  & 22.04 & \multicolumn{1}{r}{11.27} & 22.04 & 1.84  & 22.04 & \textbf{0.52} & 53.50 & 0.95  & 22.04 & 3.10  & 22.04 & 1.55  & 32.93 & 0.65  & 63.56 & 0.73  & 63.56 \\
	25    & 10    & \multicolumn{1}{r}{0.53} & 55.26 & \multicolumn{1}{r}{12.37} & 30.77 & 1.49  & 30.77 & \multicolumn{1}{r}{13.65} & 30.77 & 1.76  & 30.77 & \textbf{0.46} & 61.48 & 1.13  & 30.77 & 2.58  & 30.77 & 1.39  & 43.34 & 0.80  & 74.92 & 0.79  & 74.92 \\
	30    & 3     & \multicolumn{1}{r}{\textbf{1.03}} & 27.94 & \multicolumn{1}{r}{47.60} & 13.29 & 3.93  & 13.29 & \multicolumn{1}{r}{36.80} & 13.29 & 4.89  & 13.29 & 1.04  & 48.61 & 2.15  & 13.29 & 6.94  & 13.29 & 3.94  & 28.05 & 1.66  & 52.20 & 1.55  & 52.20 \\
	30    & 7     & \multicolumn{1}{r}{1.17} & 39.24 & \multicolumn{1}{r}{25.19} & 18.52 & 2.37  & 18.52 & \multicolumn{1}{r}{27.59} & 18.52 & 2.80  & 18.52 & \textbf{0.83} & 51.50 & 1.70  & 18.52 & 5.24  & 18.52 & 2.28  & 35.45 & 0.97  & 63.02 & 0.95  & 63.02 \\
	30    & 10    & \multicolumn{1}{r}{1.14} & 50.88 & \multicolumn{1}{r}{34.23} & 26.21 & 3.05  & 26.21 & \multicolumn{1}{r}{42.51} & 26.21 & 3.93  & 26.21 & \textbf{0.95} & 58.81 & 1.67  & 26.21 & 5.11  & 26.21 & 2.35  & 42.75 & 1.14  & 70.99 & 1.24  & 70.99 \\
	40    & 3     & \multicolumn{1}{r}{3.12} & 26.89 & \multicolumn{1}{r}{148.29} & 12.97 & 13.13 & 12.97 & \multicolumn{1}{r}{152.82} & 12.97 & 7.88  & 12.97 & 3.13  & 48.55 & 4.85  & 12.97 & 11.46 & 12.97 & 7.69  & 31.65 & 3.25  & 51.45 & \textbf{2.75} & 51.45 \\
	40    & 7     & \multicolumn{1}{r}{7.41} & 37.54 & \multicolumn{1}{r}{146.75} & 16.74 & 11.45 & 16.74 & \multicolumn{1}{r}{132.06} & 16.74 & 12.42 & 16.74 & 2.83  & 51.51 & 4.65  & 16.74 & 20.34 & 16.74 & 8.56  & 34.90 & \textbf{2.62} & 58.41 & 2.77  & 58.41 \\
	40    & 10    & \multicolumn{1}{r}{5.85} & 44.38 & \multicolumn{1}{r}{123.78} & 21.11 & 14.14 & 21.11 & \multicolumn{1}{r}{179.14} & 21.11 & 11.45 & 21.11 & 2.91  & 55.03 & 4.09  & 21.11 & 17.93 & 21.11 & 7.09  & 37.62 & \textbf{2.28} & 65.49 & 2.54  & 65.49 \\
	50    & 5     & \multicolumn{1}{r}{11.72} & 30.44 & \multicolumn{1}{r}{314.63} & 12.63 & 38.32 & 12.63 & \multicolumn{1}{r}{385.67} & 12.63 & 19.99 & 12.63 & 9.10  & 49.44 & 7.86  & 12.63 & 34.72 & 12.63 & 15.52 & 38.65 & \textbf{5.46} & 54.26 & 5.50  & 54.26 \\
	50    & 10    & \multicolumn{1}{r}{45.50} & 42.04 & \multicolumn{1}{r}{565.37} & 20.25 & 60.62 & 20.25 & \multicolumn{1}{r}{1357.87} & 20.25 & 43.47 & 20.25 & 18.39 & 54.50 & 12.96 & 20.25 & 70.69 & 20.25 & 17.71 & 43.19 & 7.79  & 62.29 & \textbf{5.94} & 62.29 \\
	50    & 15    & \multicolumn{1}{r}{57.11} & 47.35 & \multicolumn{1}{r}{294.64} & 20.88 & 32.18 & 20.88 & \multicolumn{1}{r}{615.38} & 20.88 & 41.99 & 20.88 & 9.50  & 57.07 & 10.91 & 20.88 & 36.25 & 20.88 & 15.87 & 44.26 & \textbf{6.54} & 69.01 & 8.80  & 69.01 \\
	75    & 5     & \multicolumn{1}{r}{110.47} & 28.03 & \multicolumn{1}{r}{1779.03} & 11.50 & 292.01 & 11.50 & 1304.42(1) & 11.50 & 131.39 & 11.50 & 165.08 & 48.94 & 47.01 & 11.50 & 442.62 & 11.50 & 52.17 & 48.03 & 24.48 & 54.79 & \textbf{18.73} & 54.79 \\
	75    & 10    & \multicolumn{1}{r}{927.58} & 36.27 & \multicolumn{1}{r}{2518.22} & 13.98 & 416.00 & 13.98 & 3519.38(2) & 13.98 & 173.89 & 13.98 & 239.72 & 51.55 & 100.70 & 13.98 & 239.59 & 13.98 & 66.51 & 50.36 & 26.83 & 59.35 & \textbf{19.97} & 59.35 \\
	75    & 15    & \multicolumn{1}{r}{1973.03} & 41.65 & \multicolumn{1}{r}{1931.90} & 16.57 & 380.97 & 16.57 & 3567.20(2) & 16.57 & 231.85 & 16.57 & 101.01 & 53.75 & 59.88 & 16.57 & 201.43 & 16.57 & 73.25 & 54.31 & 28.21 & 64.43 & \textbf{23.56} & 64.43 \\
	100   & 10    & 5981.32(4) & 26.05 & 3443.88(1) & 12.15 & 1090.56 & 12.15 & \multicolumn{1}{r}{(5)} & 12.15 & 476.64 & 12.15 & 457.77 & 50.15 & 164.32 & 12.15 & 605.83 & 12.15 & 186.07 & 56.74 & \textbf{63.53} & 59.06 & 64.09 & 59.06 \\
	100   & 15    & \multicolumn{1}{r}{(5)} & 37.41 & 2260.86(2) & 13.17 & 1611.93 & 13.17 & 3277.27(4) & 13.17 & 476.51 & 13.17 & 958.80 & 52.14 & 142.57 & 13.17 & 435.76 & 13.17 & 271.07 & 57.95 & 77.88 & 61.87 & \textbf{64.12} & 61.87 \\
	100   & 25    & 5098.56(4) & 49.16 & 4780.63(2) & 19.07 & 1212.76 & 19.07 & 2887.83(4) & 19.07 & 968.22 & 19.07 & 544.08 & 58.96 & 133.52 & 19.07 & 674.77 & 19.07 & 336.73 & 61.62 & 78.34 & 69.13 & \textbf{64.62} & 69.13 \\
		\hline
	\end{tabular}%
}
	\label{tab1}%
\end{sidewaystable}%

In terms of running times, observe that for $n=100$ some of the instances cannot be solved in less than two hours if formulations F1, F2 or F3 are used. However, the reported results of F2' and F3-\eqref{33}+\eqref{F22}+\eqref{F22_2} are much better than those corresponding to F2 or F3. Note that times of F4 are similar in many of the cases to those required by F3-\eqref{33}+\eqref{F22}+\eqref{F22_2} and all the instances can be solved in less than two hours.

Observe also that F5 seems to provide better results than F4. Furthermore, it is clear that the best formulation is F5 replacing constraints
\eqref{F55} by constraints \eqref{53}  and relaxing the integrality of variables $z_{ir}$ for $i\in N$, $r\in\{3,\ldots,G_i\}$. By using this variant of formulation F5, the results show that running times are (in average) not bigger than 65 seconds in any of the cases.

In contrast, the LP gaps of F2, F2', F3, F3-\eqref{33}+\eqref{F22}+\eqref{F22_2}, F5 and  F5-\eqref{F55}+\eqref{desagregada}, which always coincide, are the smallest ones. Although F5-\eqref{F55}+\eqref{53}* is the formulation that provides the best computational times, the reported LP gaps are the largest ones if we compare them with the remaining formulations.

Since  F5-\eqref{F55}+\eqref{53}* is the best formulation in terms of times, next subsection is devoted  to the computational study of this formulation reducing the number of variables and using valid inequalities.

\subsection{Reduction of variables and valid inequalities for F5-\eqref{F55}+\eqref{53}*}

In this subsection we observe the results of
using a preprocessing phase to reduce the number of variables in formulation F5-\eqref{F55}+\eqref{53}* and we will also report the results when applying valid inequalities.

In Subsection \ref{reduction} a  preprocessing phase to reduce the number of $z$- and $u$-variables is described.  Concretely, constraints \eqref{P3} allow to reduce the number of defined $z$-variables. Similarly, constraints \eqref{varured} decrease the number of $u$-variables. Besides, a reduction of $u$-variables based on obtaining an adequate lower bound of the $p$-center objective value considering each stratum independently is described.

% Table generated by Excel2LaTeX from sheet 'PREPROCESOS'
\begin{table}[!h]
	\centering
	\caption{Percentage of $z$- and $u$-variables reduced with respect  to the original ones.}
	\scalebox{0.7}{
	\begin{tabular}{rr|r|rr}
		\hline          & \multicolumn{1}{r}{} & \multicolumn{1}{|c|}{\%$z$} & \multicolumn{2}{c}{\%$u$} \\
		\hline
		\multicolumn{1}{l}{n} & \multicolumn{1}{l|}{p} & \multicolumn{1}{l|}{\eqref{P3}} & \multicolumn{1}{l}{\eqref{varured}+clas. Rel} & \multicolumn{1}{l}{\eqref{varured}+Binary alg.} \\
		\hline
		75    & 5     & 6.35  & 29.47 & \textbf{43.11} \\
		75    & 10    & 13.61 & 23.58 & \textbf{36.63} \\
		75    & 15    & 20.21 & 21.75 & \textbf{32.79} \\
		100   & 10    & 10.74 & 25.44 & \textbf{37.32} \\
		100   & 15    & 16.13 & 23.75 & \textbf{34.84} \\
		100   & 25    & 26.01 & 22.55 & \textbf{31.68} \\
		200   & 10    & 9.26  & 30.23 & \textbf{40.72} \\
		200   & 20    & 16.93 & 25.66 & \textbf{35.33} \\
		200   & 30    & 22.88 & 24.82 & \textbf{33.96} \\
		\hline
	\end{tabular}%
}
	\label{tab3}%
\end{table}%

 In particular, we mention two ways to obtain  these lower bounds. The first one
is to solve the linear relaxation  for the $pCP$ using the classic formulation of \cite{daskin1995}. The second way consists in using the binary algorithm proposed in \cite{Calik}. Table~\ref{tab3} reports the percentage of fixed $z$- and $u$-variables in formulation F5-\eqref{F55}+\eqref{53}* when the former criteria for fixing variables are applied.
The first column corresponds to the percentage of reduced $z$-variables if constraints \eqref{P3} are applied. The second column reports the percentage of fixed $u-$variables when using constraints \eqref{varured} together with the reduction strategy based on the solving of \cite{daskin1995} relaxed formulation for each stratum. Finally the last column reports the percentage of reduction when \eqref{varured} and Binary Algorithm specified in \cite{Calik} for each stratum are applied. Observe that between $6.35\%$ and $26.01\%$ of the $z$-variables could be fixed. In the case of $u$-variables the largest number of fixed $u$-variables (boldfaced) is obtained when applying the Binary Algorithm. With this strategy and \eqref{varured}, more than a $31\%$ of $u$-variables are fixed in average.

% Table generated by Excel2LaTeX from sheet 'PREPROCESOS'
\begin{table}[!ht]
	\centering
	\caption{Times and LP gaps reducing the number of $z$- and $u$-variables in formulation F5-\eqref{F55}+\eqref{53}.}
	\scalebox{0.65}{
	\begin{tabular}{rr|rr|rr|rrr|rrr|rrr}
		\hline
		&       & \multicolumn{2}{c|}{F5-\eqref{F55}+\eqref{53}} & \multicolumn{2}{c|}{F5-\eqref{F55}+\eqref{53}*} & \multicolumn{3}{c|}{Classic rel} & \multicolumn{3}{c|}{Binary} & \multicolumn{3}{c}{Binary*} \\
		\multicolumn{1}{l}{n} & \multicolumn{1}{l|}{p} & \multicolumn{1}{l}{Time} & \multicolumn{1}{l|}{LP Gap} & \multicolumn{1}{l}{Time} & \multicolumn{1}{l|}{LP Gap} & \multicolumn{1}{l}{t prepro} & \multicolumn{1}{l}{t total} & \multicolumn{1}{l|}{LP Gap} & \multicolumn{1}{l}{t prepro} & \multicolumn{1}{l}{t total} & \multicolumn{1}{l|}{LP Gap} & \multicolumn{1}{l}{t prepro} & \multicolumn{1}{l}{t total} & \multicolumn{1}{l}{LP Gap} \\
		\hline
		75    & 5     & 24.48 & 54.79 & 18.73 & 54.79 & 1.53  & 19.53 & 33.56 & 0.52  & 12.50 & \textbf{8.55} & 0.51  & \textbf{7.95} & \textbf{8.55} \\
75    & 10    & 26.83 & 59.35 & 19.97 & 59.35 & 1.09  & 16.60 & 41.45 & 0.49  & 13.56 & \textbf{10.82} & 0.49  & \textbf{11.04} & \textbf{10.82} \\
75    & 15    & 28.21 & 64.43 & 23.56 & 64.43 & 0.93  & 18.85 & 48.96 & 0.49  & 12.87 & \textbf{18.42} & 0.46  & \textbf{10.95} & \textbf{18.42} \\
100   & 10    & 63.53 & 59.06 & 64.09 & 59.06 & 2.54  & 56.19 & 39.43 & 1.08  & 42.19 & \textbf{10.94} & 1.14  & \textbf{29.22} & \textbf{10.94} \\
100   & 15    & 77.88 & 61.87 & 64.12 & 61.87 & 2.06  & 63.73 & 44.93 & 0.93  & 37.41 & \textbf{15.45} & 0.95  & \textbf{30.26} & \textbf{15.45} \\
100   & 25    & 78.34 & 69.13 & 64.62 & 69.13 & 1.63  & 43.27 & 56.49 & 0.91  & 31.84 & \textbf{23.35} & 0.91  & \textbf{28.95} & \textbf{23.35} \\
200   & 10    & 440(1) & 56.86 & 1248.75 & 56.86 & 26.75 & 739.05 & 33.83 & 9.28  & 368.87 & \textbf{8.95} & 9.19  & \textbf{275.96} & \textbf{8.95} \\
200   & 20    & 440.19 & 58.97 & 436.89 & 58.97 & 19.01 & 267.61 & 39.59 & 9.58  & 118.36 & \textbf{11.04} & 9.60  & \textbf{82.42} & \textbf{11.04} \\
200   & 30    & 349.71 & 62.75 & 503.01 & 62.75 & 13.78 & 199.57 & 46.25 & 7.80  & 111.97 & \textbf{15.28} & 7.84  & \textbf{89.68} & \textbf{15.28} \\
		\hline
	\end{tabular}%
}
	\label{tab4}%
\end{table}%
 Table \ref{tab4} reports the computational times and LP gaps  for $n\in\{75,100,200\}$ if the former  preprocessing phase for fixing variables are used in order to reduce the number of variables. The first block of columns corresponds to the formulation without any  preprocessing phase and the second one corresponds to the formulation relaxing $z_{ir}$ for $i\in N$, $r=3,\ldots, G_i$. After these two blocks, different options for the preprocessing are studied. In those cases, a first column indicating the preprocessing time is included in each block.

 Columns in block ``classic rel.'' report the results if a preprocessing using \eqref{P3} and \eqref{varured} based on the relaxed formulation from \cite{daskin1995} is used. ``Binary'' shows the results if Binary algorithm proposed in \cite{Calik} is used to obtain a lower bound on the $p$-center for each stratum and the criteria given by \eqref{P3} and \eqref{varured} are applied. In columns under heading ``Binary*'', the same  preprocessing is used but, in this case, $z_{ir}$ variables are relaxed for $i\in N$, $r=3,\ldots, G_i$.
The largest differences in CPU time among the variants can be observed in instances with $n=200$. In this case, the best results regarding CPU time are the ones reported in column ``Binary*'. It is worth noting that the preprocessing times represent only a small fraction of the overall solution time in all the instances reported in this table. Observe also that the LP gaps are considerably reduced if binary algorithm together with \eqref{P3} and \eqref{varured} is used.

% Table generated by Excel2LaTeX from sheet 'desigualdades_validas'
\begin{table}[htbp]
	\centering
	\caption{Times of F5-\eqref{F55}+\eqref{53}* using binary algorithm to reduce the number of $u$-variables and different valid inequalities.}
	\scalebox{0.65}{
	\begin{tabular}{rr|rrrrrrr}
		\hline
			\multicolumn{1}{l}{n} & \multicolumn{1}{l|}{p} & Binary*  & \eqref{R1mod}  & \eqref{escenarios}& \eqref{escenarios2}  & \eqref{escenarios3} & \eqref{Restz}  & \eqref{F2_3} \\
		\hline
		75    & 5     & 7.95  & 12.61 & 8.02  & 8.02  & 7.95  & 11.95 & \textbf{6.03} \\
		75    & 10    & 11.04 & 13.95 & 11.08 & 11.11 & \textbf{11.02} & 16.25 & 15.90 \\
		75    & 15    & 10.95 & 12.96 & 11.09 & 11.08 & 11.01 & 14.24 & \textbf{9.44} \\
		100   & 10    & 29.22 & 37.42 & 29.21 & 29.13 & 29.27 & 43.85 & \textbf{24.17} \\
		100   & 15    & 30.26 & 45.75 & 30.29 & 30.28 & 30.22 & 46.58 & \textbf{28.28} \\
		100   & 25    & 28.95 & 37.38 & 28.82 & \textbf{28.73} & 28.76 & 44.35 & 36.95 \\
		200   & 10    & 275.96 & 539.47 & 275.58 & 275.37 & 276.19 & 289.21 & \textbf{162.02} \\
		200   & 20    & \textbf{82.42} & 161.46 & 82.49 & 82.69 & 82.66 & 92.31 & 93.18 \\
		200   & 30    & \textbf{89.68} & 176.37 & 90.28 & 90.04 & 90.00 & 120.44 & 164.50 \\
		300   & 15    & 509.79 & 1298.54 & 512.82 & 513.38 & 510.10 & 523.33 & \textbf{271.05} \\
		300   & 30    & 315.13 & 591.42 & 318.61 & 316.23 & 316.46 & 372.64 & \textbf{228.18} \\
		300   & 45    & 535.69 & 813.88 & 538.52 & 533.47 & 532.62 & \textbf{442.23} & 610.77 \\
		400   & 20    & 1017.28 & 3305.29 & 1011.01 & 1012.90 & 1014.40 & 722.30 & \textbf{450.12} \\
		400   & 40    & 663.16 & 1863.28 & 666.45 & \textbf{660.53} & 663.98 & 805.02 & 954.81 \\
		400   & 60    & 475.14 & 1246.22 & \textbf{474.36} & 475.05 & 474.23 & 735.84 & 816.77 \\
		\hline
	\end{tabular}%
}
	\label{tab2}%
\end{table}%
Table \ref{tab2} reports the average times required to solve the same instances with formulation  F5-\eqref{F55}+\eqref{53}*  using Binary Algorithm, \eqref{P3} and \eqref{varured} to reduce the number of variables  and  adding some of the constraints explained in Subsection  \ref{inequalitiessection}.
Regarding the reported results in Table \ref{tab2}, the time performance is significantly improved in some cases if constraints \eqref{F2_3} are included as valid inequalities for the formulation. The remaining valid inequalities appearing in this table, except maybe for \eqref{R1mod}, do not worsen the times in general, but they neither provide a significant improvement.

Finally, Table \ref{tab7} reports the time results using ORLIB data with the same $p$ values as in the original instances and using random strata. For solving these instances, formulation F5-\eqref{F55}+\eqref{53}* was used with Binary Algorithm and adding  \eqref{P3} and \eqref{varured} to reduce the number of variables. The results shows that only two instances remain unsolved after two hours using the model with the proposed preprocessing phase (underlined cpu time). In this table, we give separately the time to solve the formulation, under heading $t_{\textrm{solv}}$, the preprocessing time, under $t_{\textrm{prep}}$ and the overall time, $t_{\textrm{total}}$. Additionally, we provide the number of nodes explored in the branch and bound tree.}

% Table generated by Excel2LaTeX from sheet 'ORLIB'
\begin{sidewaystable}[!h]
	\centering
	\caption{Results for ORLIB data}
	\scalebox{0.7}{
		\begin{tabular}{rrr|rrrrr||rrr|rrrrr}
			\hline
			& \multicolumn{1}{l}{n} & \multicolumn{1}{l|}{p} & \multicolumn{1}{l}{$t_{\mbox{solv}}$} & \multicolumn{1}{l}{$t_{\mbox{prep}}$} & \multicolumn{1}{l}{$t_{\mbox{total}}$} & \multicolumn{1}{l}{LP Gap} & \multicolumn{1}{l||}{\# Nodes} &       & \multicolumn{1}{l}{n} & \multicolumn{1}{l|}{p} & \multicolumn{1}{l}{$t_{\mbox{solv}}$} & \multicolumn{1}{l}{$t_{\mbox{prep}}$} & \multicolumn{1}{l}{$t_{\mbox{total}}$} & \multicolumn{1}{l}{LP Gap} & \multicolumn{1}{l}{\# Nodes} \\
			\hline
 pmed1  & 100 &   5 &   46.00 &   1.12 &   47.66 &  8.75 &  8041 & pmed21 & 500 &   5 &  120.12 & 334.01 & 466.05 & 4.86  & 1387 \\
 pmed2  & 100 &  10 &   26.08 &   1.14 &   28.13 &  9.91 &  4783 & pmed22 & 500 &  10 &  253.44 & 307.79 & 597.42 & 7.78  & 1713 \\
 pmed3  & 100 &  10 &   11.51 &   1.09 &   12.98 &  8.05 &  1271 & pmed23 & 500 &  50 & \underline{7200.19} & 229.07 & 7447.45 & 18.36 & 29644 \\
 pmed4  & 100 &  20 &   28.27 &   0.92 &   29.48 & 18.50 &  9463 & pmed24 & 500 & 100 & 1068.26 & 320.28 & 1392.83 & 22.16 & 18441 \\
 pmed5  & 100 &  33 &   18.61 &   0.72 &   19.52 & 37.80 &  5843 & pmed25 & 500 & 167 & 1591.44 & 199.79 & 1793.49 & 40.16 & 81497 \\
 pmed6  & 200 &   5 &   75.99 &  10.02 &   92.17 &  6.95 &  2621 & pmed26 & 600 &   5 &  259.63 & 677.29 & 959.90 & 4.58  & 1331 \\
 pmed7  & 200 &  10 &   89.48 &   8.24 &   99.64 &  8.59 &  4737 & pmed27 & 600 &  10 &  427.28 & 593.22 & 1056.32 & 8.91  & 2315 \\
 pmed8  & 200 &  20 &  279.66 &   8.63 &  290.36 & 13.40 & 10471 & pmed28 & 600 &  60 & 1784.88 & 572.43 & 2371.33 & 13.53 & 6377 \\
 pmed9  & 200 &  40 &   48.49 &   8.75 &   57.95 & 17.09 &  4737 & pmed29 & 600 & 120 &  487.86 & 425.46 & 918.56 & 22.71 & 10189 \\
 pmed10 & 200 &  67 &   47.67 &   6.26 &   54.21 & 33.24 & 10931 & pmed30 & 600 & 200 &  188.01 & 505.15 & 696.21 & 32.72 & 5877 \\
 pmed11 & 300 &   5 &   30.98 &  44.26 &   78.77 &  5.05 &   353 & pmed31 & 700 &   5 &  132.35 & 1318.46 & 1470.65 & 2.73  & 393 \\
 pmed12 & 300 &  10 &   57.64 &  41.64 &  105.20 &  6.77 &  1215 & pmed32 & 700 &  10 & 4250.54 & 788.62 & 5104.98 & 7.09  & 19239 \\
 pmed13 & 300 &  30 &  786.58 &  33.21 &  829.07 & 17.75 & 13347 & pmed33 & 700 &  70 & 2666.15 & 1024.91 & 3729.87 & 13.66 & 7367 \\
 pmed14 & 300 &  60 &  338.51 &  38.84 &  380.14 & 19.62 &  9919 & pmed34 & 700 & 140 &  305.98 & 752.26 & 1067.68 & 16.55 & 3559 \\
 pmed15 & 300 & 100 &   60.15 &  40.86 &  101.79 & 27.72 &  6159 & pmed35 & 800 &   5 &  180.27 & 1429.91 & 1647.11 & 4.18  & 935 \\
 pmed16 & 400 &   5 &   47.96 & 114.38 &  165.70 &  3.88 &   599 & pmed36 & 800 &  10 & 2376.87 & 1837.73 & 4283.38 & 6.03  & 4411 \\
 pmed17 & 400 &  10 &  184.14 & 115.38 &  309.41 &  6.19 &  2399 & pmed37 & 800 &  80 & 5553.53 & 1353.61 & 6987.60 & 14.52 & 9285 \\
 pmed18 & 400 &  40 &  384.20 &  86.96 &  480.31 &  9.80 &  4275 & pmed38 & 900 &   5 &  207.01 & 2844.10 & 3109.63 & 5.43  & 289 \\
 pmed19 & 400 &  80 & 3354.26 & 119.21 & 3476.32 & 21.85 & 64273 & pmed39 & 900 &  10 & 2993.05 & 2981.84 & 6021.92 & 7.82  & 6587 \\
 pmed20 & 400 & 133 &  137.76 &  93.43 &  232.63 & 39.05 & 10437 & pmed40 & 900 &  90 & \underline{7205.28} & 3780.23 & 11033.40 & 13.07 & 10222 \\
\hline
\end{tabular}%
	}
	\label{tab7}%
\end{sidewaystable}%

In this table we observe that varying $p$ has a strong effect on the CPU times, both, in the preprocessing phase and when solving the final formulation. Moreover, the effect is different in both cases, yielding curious situations, where the preprocessing time can be larger than the actual solution time.
\begin{figure}[!h]
\begin{center}
\includegraphics[width=.45\textwidth]{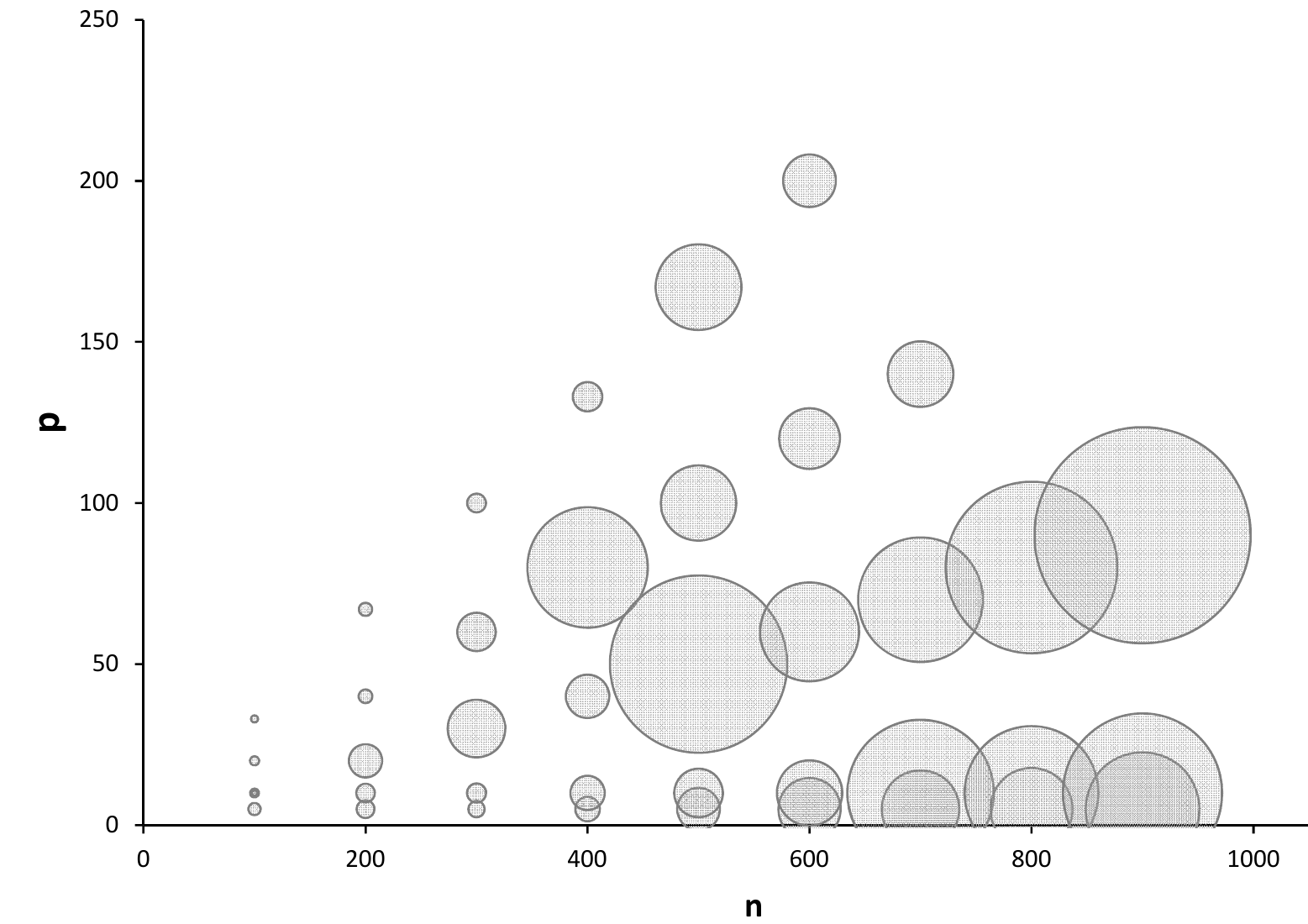}
\caption{CPU times (circle size) as a function of $n$ and $p$ \label{fig:ltimes}}
\end{center}
\end{figure}
We can also observe that the most demanding instances tend to be those wit $p\simeq 10\%\cdot n$. This behavior can be better appreciated in Figure~\ref{fig:ltimes}.
\subsection{SAA for P$p$CP}
 In this subsection, the time and gap results of SAA for the P$p$CP are analyzed. Table \ref{tab5} shows the results of SAA in comparison with P$p$CP formulation presented in \cite{PpCP}.

The first column corresponds to the running time of the probability chain P$p$CP formulation described in \cite{PpCP}. ``F1 SAA'' shows the results of SAA if formulation F1 of the S$p$CP is used.  ``Binary* SAA'' resports again the results of SAA but using formulation F5 with constraints \eqref{F55} replaced by \eqref{53}, using Binary Algorithm as a preprocessing phase and relaxing $z_{ir}$ variables for $i\in N$ and $r\in\{3,\ldots,G_i\}$ .   For each block of columns, the gap column reports the gap (in percetage, $\%$) between  the best obtained solution in the SAA heuristic and the P$p$CP objective value. In addition, the time column reports the running time of the procedures.

 Regarding the running times of SAA, we observe a significant difference between SAA when using formulation F1 and the remaining SAA columns that use formulation F5. As observed,  times in ``Binary* SAA'' grow much slower than when using F1 so that, even if for the smallest instances they seem to be worse, they become much better for $n > 30$.
Considering the gaps we see that in none of the cases, the gaps are bigger than $0.64\%$. Moreover, both versions of the SAA found the optimal solution for at least half of the instances.  As explained in Section \ref{saa}, we can find theoretical results that guarantee the goodness of the obtained solution when using the SAA.

 Table \ref{tab6} reports the average results of the instances with $(n,p)\in\{(50,15),(75,10),(100,10),$ $(100,15),(100,25)\}$. First column reports the necessary time for solving the P$p$CP using the probability chain formulation.
	``Gap$_{BS}$'' column reports the gap between the best solution obtained by SAA method and the best solution of P$p$CP within the time limit. Finally, SAA time is reported. Observe that in all unsolved instances after two hours ``Gap$_{BS}$'' column reports negative gaps. This is due to the fact that the best solution given by SAA is better than the best solution provided by P$p$CP formulation after two hours.

% Table generated by Excel2LaTeX from sheet 'Hoja1'
\begin{table}[htbp]
	\centering
	\caption{SAA results}
	\scalebox{0.7}{
	\begin{tabular}{rr|r|rr|rr}
		\hline
		&       & \multicolumn{1}{l|}{P$p$CP} & \multicolumn{2}{c|}{F1 SAA} & \multicolumn{2}{c}{Binary* SAA} \\
		\multicolumn{1}{r}{n} & \multicolumn{1}{r|}{p} & \multicolumn{1}{r|}{Time} & \multicolumn{1}{r}{Gap} & \multicolumn{1}{r|}{Time} & \multicolumn{1}{r}{Gap} & \multicolumn{1}{r}{Time} \\
		\hline
		6     & 2     & \textbf{0.01} & 0.00  & 0.10  & 0.00  & 0.45 \\
		10    & 3     & \textbf{0.03} & 0.00  & 0.54  & 0.00  & 0.87 \\
		10    & 5     & \textbf{0.03} & 0.00  & 0.49  & 0.03  & 1.49 \\
		13    & 3     & \textbf{0.07} & 0.00  & 0.91  & 0.00  & 1.53 \\
		13    & 5     & \textbf{0.07} & 0.00  & 1.10  & 0.24  & 2.21 \\
		13    & 8     & \textbf{0.05} & 0.00  & 1.32  & 0.00  & 2.20 \\
		15    & 3     & \textbf{0.10} & 0.00  & 1.77  & 0.00  & 1.88 \\
		15    & 7     & \textbf{0.13} & 0.00  & 1.73  & 0.10  & 3.56 \\
		15    & 10    & \textbf{0.07} & 0.00  & 1.84  & 0.56  & 2.65 \\
		20    & 3     & \textbf{0.32} & 0.00  & 3.40  & 0.00  & 4.18 \\
		20    & 7     & \textbf{0.63} & 0.64  & 3.67  & 0.00  & 8.49 \\
		20    & 10    & \textbf{0.49} & 0.10  & 4.28  & 0.14  & 5.62 \\
		25    & 3     & \textbf{0.84} & 0.00  & 6.73  & 0.00  & 7.42 \\
		25    & 7     & \textbf{3.48} & 0.05  & 8.96  & 0.23  & 9.50 \\
		25    & 10    & \textbf{5.13} & 0.02  & 9.48  & 0.01  & 14.34 \\
		30    & 3     & \textbf{2.01} & 0.00  & 13.90 & 0.00  & 11.27 \\
		30    & 7     & 13.61 & 0.14  & 12.78 & 0.15  & \textbf{9.40} \\
		30    & 10    & 22.99 & 0.00  & \textbf{16.24} & 0.00  & 16.54 \\
		40    & 3     & \textbf{8.28} & 0.00  & 40.90 & 0.00  & 19.94 \\
		40    & 7     & 148.22 & 0.01  & 98.39 & 0.20  & \textbf{19.45} \\
		40    & 10    & 295.52 & 0.01  & 96.68 & 0.01  & \textbf{19.52} \\
		50    & 5     & 243.17 & 0.03  & 162.76 & 0.00  & \textbf{44.68} \\
		50    & 10    & 4083.75 & 0.01  & 462.26 & 0.12  & \textbf{67.74} \\
		75    & 5     & 4108.22 & 0.03  & 1386.77 & 0.03  & \textbf{150.28} \\
		\hline
	\end{tabular}%
}
	\label{tab5}
\end{table}%

\clearpage

% Table generated by Excel2LaTeX from sheet 'PRUEBA3'
\begin{table}[htbp]
	\centering
	\caption{SAA results for larger instances.}
	\scalebox{0.7}{
	\begin{tabular}{rr|lrr}
		\hline
		\multicolumn{1}{l}{n} & \multicolumn{1}{l|}{p} & P$p$CP Time & \multicolumn{1}{l}{Gap$_{BS}$} & \multicolumn{1}{l}{SAA Time} \\
		\hline
		50    & 15    & 3781.82(3) & -1.52 & 71.63 \\
		75    & 10    & $>$7200 & -9.55 & 200.57 \\
		75    & 15    & $>$7200 & -13.03 & 258.32 \\
		100   & 10    & $>$7200 & -12.98 & 491.73 \\
		100   & 15    & $>$7200 & -17.60 & 449.02 \\
		100   & 25    & $>$7200 & -21.71 & 850.55 \\
		\hline
	\end{tabular}%
}
	\label{tab6}%
\end{table}%

\section{Conclusions}\label{conclusion}

	This paper presents an extension of the $p$-center problem
	called the Stratified $p$-Center Problem (\spcp). This extension could be
	applied in cases where the population is divided into different strata and
	the evaluation of the service must be separately measured for each stratum.
	In the model, it is assumed that more than one stratum can be
	present at each demand point. 
	
	Different formulations were introduced together with a detailed study of variants, variable reduction
	processes and valid inequalities. Regarding the computational results, the
	best performance was obtained using a formulation based on covering
	variables.
	
	The \spcp allows to implement a heuristic approach based on the
	Sample Average Approximation (SAA) method to obtain good feasible
	solutions for the probabilistic $p$-center problem. This heuristic approach provides good upper bounds in acceptable times.
	
\section*{Acknowledgements}
A.M. Rodr{\'\i}guez-Ch{\'\i}a and Luisa I. Mart{\'\i}nez Merino   acknowledge that research reported here was supported by  the European Regional Development's funds (FEDER) and Agencia Estatal de Investigaci\'on (AEI)  under project MTM2016-74983-C2-2-R. Luisa I. Mart\'inez Merino was also supported by Universidad de C\'adiz PhD grant UCA/REC02VIT/2014 and Programa de Fomento e Impulso de la actividad Investigadora UCA (2018). The research of Maria Albareda has been partially funded by Spanish Ministry of Economy and Competitiveness end EDRF funds through project MTM2015-63779-R.

\bibliographystyle{abbrvnat}
%\bibliography{references}

\begin{thebibliography}{25}
	\providecommand{\natexlab}[1]{#1}
	\providecommand{\url}[1]{\texttt{#1}}
	\expandafter\ifx\csname urlstyle\endcsname\relax
	\providecommand{\doi}[1]{doi: #1}\else
	\providecommand{\doi}{doi: \begingroup \urlstyle{rm}\Url}\fi
	
	\bibitem[Albareda-Sambola et~al.(2010)Albareda-Sambola, D{\'\i}az, and
	Fern{\'a}ndez]{capacPCP2}
	M.~Albareda-Sambola, J.~D{\'\i}az, and E.~Fern{\'a}ndez.
	\newblock Lagrangean duals and exact solution to the capacitated p-center
	problem.
	\newblock \emph{European Journal of Operational Research}, 201\penalty0
	(1):\penalty0 71--81, 2010.
	
	\bibitem[Averbakh and Berman(1997)]{robust2}
	I.~Averbakh and O.~Berman.
	\newblock Minimax regret p-center location on a network with demand
	uncertainty.
	\newblock \emph{Location Science}, 5\penalty0 (4):\penalty0 247--254, 1997.
	
	\bibitem[Balinski(1965)]{MIP1}
	M.~Balinski.
	\newblock Integer programming: methods uses, computation.
	\newblock \emph{Management Science}, 12:\penalty0 253--313, 1965.
	
	\bibitem[Calik and Tansel(2013)]{Calik}
	H.~Calik and B.~C. Tansel.
	\newblock Double bound method for solving the p-center location problem.
	\newblock \emph{Computers and Operations Research}, 40\penalty0 (12):\penalty0
	2991--2999, 2013.
	
	\bibitem[Calik et~al.(2015)Calik, Labb{\'e}, and Yaman]{calik2015p}
	H.~Calik, M.~Labb{\'e}, and H.~Yaman.
	\newblock p-center problems.
	\newblock In \emph{Location Science}, pages 79--92. Springer, 2015.
	\newblock ISBN 978-3-319-13110-8.
	
	\bibitem[Daskin(1995)]{daskin1995}
	M.~Daskin.
	\newblock \emph{Network and Discrete Location: Models, Algorithms, and
		Applications}.
	\newblock Wiley, New York, 1995.
	
	\bibitem[Dobson and Karmarkar(1987)]{Dobson}
	G.~Dobson and U.~Karmarkar.
	\newblock Competitive location on a network.
	\newblock \emph{Operational Research}, 35:\penalty0 565--574, 1987.
	
	\bibitem[Drezner(1989)]{conditional}
	Z.~Drezner.
	\newblock Conditional p-center problems.
	\newblock \emph{Transportation Science}, 23\penalty0 (1):\penalty0 51--53,
	1989.
	
	\bibitem[Elloumi et~al.(2004)Elloumi, Labb\'{e}, and Pochet]{Elloumi}
	S.~Elloumi, M.~Labb\'{e}, and Y.~Pochet.
	\newblock A new formulation and resolution method for the p-center problem.
	\newblock \emph{INFORMS Journal on Computing}, 16\penalty0 (1):\penalty0
	84--94, 2004.
	
	\bibitem[Espejo et~al.(2012)Espejo, Mar\'{\i}n, and
	Rodr\'{\i}guez-Ch\'{\i}a]{cac}
	I.~Espejo, A.~Mar\'{\i}n, and A.~M. Rodr\'{\i}guez-Ch\'{\i}a.
	\newblock Closest assignment constraints in discrete location problems.
	\newblock \emph{European Journal Operational Research}, 219:\penalty0 49--58,
	2012.
	
	\bibitem[Espejo et~al.(2015)Espejo, Mar{\'\i}n, and
	Rodr{\'\i}guez-Ch{\'\i}a]{StocCap}
	I.~Espejo, A.~Mar{\'\i}n, and A.~M. Rodr{\'\i}guez-Ch{\'\i}a.
	\newblock Capacitated p-center problem with failure foresight.
	\newblock \emph{European Journal of Operational Research}, 247\penalty0
	(1):\penalty0 229--244, 2015.
	
	\bibitem[Garc\'ia et~al.(2011)Garc\'ia, Labb\'e, and Mar\'in]{Garcia}
	S.~Garc\'ia, M.~Labb\'e, and A.~Mar\'in.
	\newblock Solving large p-median problems with a radius formulation.
	\newblock \emph{INFORMS Journal on Computing}, 23\penalty0 (4):\penalty0
	546--556, 2011.
	
	\bibitem[Garfinkel et~al.(1977)Garfinkel, Neebe, and Rao]{garfinkelPC1}
	R.~Garfinkel, A.~Neebe, and M.~Rao.
	\newblock The m-center problem: Minimax facility location.
	\newblock \emph{Management Science}, 23\penalty0 (10):\penalty0 1133--1142,
	1977.
	
	\bibitem[{Homem-de-Mello} and Bayraksan(2014)]{Homem}
	T.~{Homem-de-Mello} and G.~Bayraksan.
	\newblock Monte carlo sampling-based methods for stochastic optimization.
	\newblock \emph{Surveys in Operations Research and Management Science},
	19\penalty0 (1):\penalty0 56--85, 2014.
	
	\bibitem[Kleywegt et~al.(2002)Kleywegt, Shapiro, and {Homem-de
		Mello}]{Kleiwegt}
	A.~Kleywegt, A.~Shapiro, and T.~{Homem-de Mello}.
	\newblock The sample average approximation method for stochastic discrete
	optimization.
	\newblock \emph{SIAM Journal on Optimization}, 12\penalty0 (2):\penalty0
	479--502, 2002.
	
	\bibitem[Linderoth et~al.(2006)Linderoth, Shapiro, and Wright]{Linderoth}
	J.~Linderoth, A.~Shapiro, and S.~Wright.
	\newblock The empirical behavior of sampling methods for stochastic
	programming.
	\newblock \emph{Annals of Operations Research}, 142\penalty0 (1):\penalty0
	215--241, 2006.
	
	\bibitem[Lu and Sheu(2013)]{robust}
	C.-C. Lu and J.-B. Sheu.
	\newblock Robust vertex p-center model for locating urgent relief distribution
	centers.
	\newblock \emph{Computers \& Operations Research}, 40\penalty0 (8):\penalty0
	2128--2137, 2013.
	
	\bibitem[Mar\'in et~al.(2009)Mar\'in, Nickel, Puerto, and Velten]{Marin}
	A.~Mar\'in, S.~Nickel, J.~Puerto, and S.~Velten.
	\newblock A flexible model and efficient solution strategies for discrete
	location problems.
	\newblock \emph{Discrete Applied Mathematics}, 157\penalty0 (5):\penalty0
	1128--1145, 2009.
	
	\bibitem[Mart{\'\i}nez-Merino et~al.(2017)Mart{\'\i}nez-Merino,
	Albareda-Sambola, and Rodr{\'\i}guez-Ch{\'\i}a]{PpCP}
	L.~I. Mart{\'\i}nez-Merino, M.~Albareda-Sambola, and A.~M.
	Rodr{\'\i}guez-Ch{\'\i}a.
	\newblock The probabilistic p-center problem: Planning service for potential
	customers.
	\newblock \emph{European Journal of Operational Research}, 262\penalty0
	(2):\penalty0 509--520, 2017.
	
	\bibitem[{\"O}zsoy and P{\i}nar(2006)]{capacPCP}
	F.~A. {\"O}zsoy and M.~{\c{C}}. P{\i}nar.
	\newblock An exact algorithm for the capacitated vertex p-center problem.
	\newblock \emph{Computers \& Operations Research}, 33\penalty0 (5):\penalty0
	1420--1436, 2006.
	
	\bibitem[Revelle and Hogan(1989)]{SPpcenter}
	C.~Revelle and K.~Hogan.
	\newblock The maximum reliability location problem and
	$\alpha$-reliablep-center problem: Derivatives of the probabilistic location
	set covering problem.
	\newblock \emph{Annals of Operations Research}, 18\penalty0 (1):\penalty0
	155--173, 1989.
	
	\bibitem[Robinson(1996)]{Robinson}
	S.~M. Robinson.
	\newblock Analysis of sample-path optimization.
	\newblock \emph{Mathematics of Operations Research}, 21\penalty0 (3):\penalty0
	513--528, 1996.
	
	\bibitem[Rubinstein and Shapiro(1990)]{Rubinstein}
	R.~Y. Rubinstein and A.~Shapiro.
	\newblock Optimization of static simulation models by the score function
	method.
	\newblock \emph{Mathematics and Computers in Simulation}, 32\penalty0
	(4):\penalty0 373--392, 1990.
	
	\bibitem[Schilling et~al.(1979)Schilling, Elzinga, Cohon, Church, and
	ReVelle]{schilling1979}
	D.~Schilling, D.~J. Elzinga, J.~Cohon, R.~Church, and C.~ReVelle.
	\newblock The team/fleet models for simultaneous facility and equipment siting.
	\newblock \emph{Transportation Science}, 13\penalty0 (2):\penalty0 163--175,
	1979.
	
	\bibitem[Shapiro(2013)]{Shapiro}
	A.~Shapiro.
	\newblock Sample average approximation.
	\newblock In \emph{Encyclopedia of Operations Research and Management Science},
	pages 1350--1355. Springer, 2013.
	
\end{thebibliography}

\end{document}